\documentclass[a4paper, 12pt]{article}
\usepackage{mathtools}
\usepackage{multirow}
\usepackage{bm}
\usepackage{amssymb}
\usepackage{amsthm}
\usepackage{url}
\usepackage{ulem}
\usepackage{graphics}

\theoremstyle{thmstyleone}
\newtheorem{theorem}{Theorem}

\theoremstyle{remark}

\newtheorem{step}{Step}

\DeclareMathOperator{\argmin}{argmin}
\begin{document}
\title{Hypothesis testing for varying coefficient models in tail index regression}
\author{{\sc{Koki Momoki}}$^{1}$\thanks{E-mail: k3499390@kadai.jp} and {\sc{Takuma Yoshida}}$^{2}$\\\\
$^{1,2}${\it{Graduate School of Science and Engineering, Kagoshima University}}\\
{\it{1-21-40 Korimoto, Kagoshima, Kagoshima, 890-8580, Japan}}}

\date{}

\maketitle

\abstract{This study examines the varying coefficient model in tail index regression. The varying coefficient model is an efficient semiparametric model that avoids the curse of dimensionality when including large covariates in the model. In fact, the varying coefficient model is useful in mean, quantile, and other regressions. The tail index regression is not an exception. However, the varying coefficient model is flexible, but leaner and simpler models are preferred for applications. Therefore, it is important to evaluate whether the estimated coefficient function varies significantly with covariates. If the effect of the non-linearity of the model is weak, the varying coefficient structure is reduced to a simpler model, such as a constant or zero. Accordingly, the hypothesis test for model assessment in the varying coefficient model has been discussed in mean and quantile regression. However, there are no results in tail index regression. In this study, we investigate the asymptotic properties of an estimator and provide a hypothesis testing method for varying coefficient models for tail index regression.}

\bigskip

\noindent
{\it{Extreme value theory, Hypothesis testing, Pareto-type model, Tail index regression, Varying coefficient model}}

\section{Introduction}\label{sec1}

In various fields, predicting the high- or low-tail behavior of data distribution is of interest. Examples include events such as heavy rains, large earthquakes, and significant fluctuations in stock prices. Extreme value theory is a standard approach for analyzing the data of such extremal events. Let $Y_1, Y_2, \ldots, Y_n$ be independent and identically distributed random variables with distribution function $F$. In extreme value theory, the following maximum domain of attraction assumption is standard: Assume that there exist sequences of constants $a_n>0$ and $b_n\in\mathbb{R}$ and a non-degenerate distribution function $G$ such that
\begin{equation*}
	\lim_{n\to\infty}P\left(\frac{\max(Y_1, \ldots, Y_n)-b_n}{a_n}\leq y\right)=G(y)
\end{equation*}
for each continuity point $y$ in $G$. This assumption implies that there exist a constant $\gamma\in\mathbb{R}$ and a positive function $\sigma(t)$ such that
\begin{equation*}
	\lim_{t\uparrow y^*}P\left(\frac{Y_i-t}{\sigma(t)}>y\mid Y_i>t\right)=(1+\gamma y)^{-1/\gamma}
\end{equation*}
for all $y$ for which $1+\gamma y>0$, where $y^*=\sup\{y: F(y)<1\}\in(-\infty, \infty]$ and the right-hand side for $\gamma=0$ is interpreted as $e^{-y}$ (see, Theorem 1.1.6 of de Haan and Ferreira 2006). The class of distributions on the right-hand side is called the generalized Pareto distribution and the parameter $\gamma$ is called the extreme value index. Therefore, in extreme value theory, the tail behavior of the data distribution is characterized by the extreme value index $\gamma$. Its existing estimators include the Hill estimator (Hill 1975), Pickands estimator (Pickands 1975), kernel estimator (Csorgo et al. 1985), maximum likelihood estimator (Smith 1987), and moment estimator (Dekkers et al. 1989), etc. It is noteworthy that the generalized Pareto distribution has different features depending on the sign of $\gamma$. If $\gamma>0$, we have
\begin{equation*}
	1-F(y)=y^{-1/\gamma}\mathcal{L}(y)
\end{equation*}
for all $y>0$, where $\mathcal{L}(y)$ is a slowly varying function at infinity; i.e., $\mathcal{L}(ys)/\mathcal{L}(y)\to1$ as $y\to\infty$ for all $s>0$. The class of these distributions is called the Pareto-type distribution. This case seems to be common in areas such as finance and insurance, and we frequently observe extremely large values in the data compared to the case of $\gamma\leq0$. Therefore, many researchers in extreme value theory have focused on this case. The Hill estimator mentioned above is one of the estimators of the positive extreme value index $\gamma$ and is widely used in many extreme value studies. In this study, we assume that the extreme value index $\gamma$ is positive.

In recent years, various regression models of the conditional extreme value index have been studied in the so-called tail index regression, where the tail index refers to the inverse of the extreme value index. The nonparametric tail index regression estimators include Gardes and Girard (2010), Stupfler (2013), Daouia et al. (2013), Gardes and Stupfler (2014), Goegebeur et al. (2014, 2015), and Ma et al. (2020). For fully nonparametric methods, the curse of dimensionality arises when multiple covariates are used. However, in many applications, extremal events are often triggered by multiple factors, and we hope to consider these factors. To avoid the curse of dimensionality, Wang and Tsai (2009) studied the parametric tail index regression assuming the linear model. However, in some applications of extreme value theory, the linear model may be too simple to predict the tail behavior of the distribution of the response. As an extension of Wang and Tsai (2009), Youngman (2019) studied additive models, Li et al. (2022) developed partially linear models, Yoshida (2022) provided single-index models, and Ma et al. (2019) provided varying coefficient models. The varying coefficient model is useful for analyzing time series and longitudinal data, etc. Because time or location is often important in many applications of extreme value theory, the varying coefficient model is expected to be useful in tail index regression. We are also interested in tail index regression assuming the varying coefficient model. 

The varying coefficient models pioneered by Hastie and Tibshirani (1993) assume that the regression function $m_Y({\bf{X}}, {\bf{T}})$ of interest satisfies
\begin{equation*}
	m_Y({\bf{X}}, {\bf{T}})={\bf{X}}^\top{\bm{\theta}}({\bf{T}})
\end{equation*}
for the given explanatory variable vectors ${\bf{X}}$ and ${\bf{T}}$, and the response variable $Y$, where ${\bm{\theta}}(\cdot)$ is the vector of unknown smooth functions of ${\bf{T}}$, which is denoted by the coefficient function vector. In mean and quantile regression, many authors have developed varying coefficient models, such as those of Wu et al. (1998), Fan and Zhang (1999), Huang et al. (2002, 2004), Kim (2007), Cai and Xu (2008), and Andriyana et al. (2014, 2018). Fan and Zhang (2008) provided a review article on the varying coefficient model. Some of these studies examined not only the estimation methods of the coefficient function, but also the hypothesis testing methods. We denote ${\bm{\theta}}(\cdot)=(\theta_1(\cdot), \theta_2(\cdot), \ldots, \theta_p(\cdot))^\top$. The hypothesis test for the constancy of a specific component can be represented as
\begin{equation*}
	{\rm{H}}_{0{\rm{C}}}: \theta_j(\cdot)\equiv C_0\quad{\rm{vs.}}\quad{\rm{H}}_{1{\rm{C}}}: \theta_j(\cdot)\not\equiv C_0
\end{equation*}
for an unknown constant $C_0$, where ${\rm{H}}_{0{\rm{C}}}$ is the null hypothesis and ${\rm{H}}_{1{\rm{C}}}$ is the alternative hypothesis. It is particularly important to test the sparsity of a specific covariate, which can be expressed as
\begin{equation*}
	{\rm{H}}_{0{\rm{Z}}}: \theta_j(\cdot)\equiv0\quad{\rm{vs.}}\quad{\rm{H}}_{1{\rm{Z}}}: \theta_j(\cdot)\not\equiv0,
\end{equation*}
where ${\rm{H}}_{0{\rm{Z}}}$ is the null hypothesis and ${\rm{H}}_{1{\rm{Z}}}$ is the alternative hypothesis. The varying coefficient model is flexible, but simpler models provide an easy interpretation of the data structure in real data analysis. The above hypothesis tests help to reduce the varying coefficient model to a simpler model. In mean and quantile regression, testing methods based on the comparison of the residual sum of squares include Cai et al. (2000), Fan et al. (2001), Huang et al. (2002), and Kim (2007), among others, where they used the bootstrap to implement their test. In mean regression, Fan and Zhang (2000) proposed the testing method based on the asymptotic distribution of the maximum deviation between the estimated coefficient function and true coefficient function.

In this study, we employ a logarithmic transformation to link the extreme value index of the response variable $Y$ to the explanatory variable vectors ${\bf{X}}$ and ${\bf{T}}$ via
\begin{equation*}
	\log\gamma({\bf{X}}, {\bf{T}})^{-1}={\bf{X}}^\top{\bm{\theta}}({\bf{T}}).
\end{equation*}
To the best of our knowledge, Ma et al. (2019) also studied this model. They provided a kernel-type nonparametric estimator of ${\bm{\theta}}({\bf{T}})$ and established asymptotic normality. However, they did not discuss hypothesis testing. Therefore, there are no results for the hypothesis tests in tail index regression. Our study aims to establish a testing method for varying coefficient models for tail index regression.

The remainder of this paper is organized as follows. Section \ref{sec2} introduces the local constant (Nadaraya-Watson type) maximum likelihood estimator of the coefficient functions, and Section \ref{sec3} investigates its asymptotic properties. Section \ref{sec4} introduces the proposed method for testing the structure of the coefficient functions and demonstrates the finite sample performance through simulations. A real example is analyzed in Section \ref{sec5}. All technical details are provided in Appendix.

\section{Model and method}\label{sec2}

\subsection{Varying coefficient models in tail index regression}\label{subsec2.1}

Let $Y>0$ be the univariate response variable of interest, ${\bf{X}}=(X_1, X_2, \ldots, X_p)^\top\in\mathbb{R}^p$ be the $p$-dimensional explanatory variable vector, and ${\bf{T}}=(T_1, T_2, \ldots, T_q)^\top\in\mathbb{R}^q$ be the $q$-dimensional explanatory variable vector. In addition, let $F(y\mid {\bf{x}}, {\bf{t}})=P(Y\leq y\mid {\bf{X}}={\bf{x}}, {\bf{T}}={\bf{t}})$ be the conditional distribution function of $Y$ given $({\bf{X}}, {\bf{T}})=({\bf{x}}, {\bf{t}})$. We consider the Pareto-type distribution
\begin{equation}
	1-F(y\mid {\bf{x}}, {\bf{t}})=y^{-1/\gamma({\bf{x}}, {\bf{t}})}\mathcal{L}(y; {\bf{x}}, {\bf{t}}),\label{2.1.1}\tag{2.1}
\end{equation}
where $\gamma({\bf{x}}, {\bf{t}})$ is a positive function of ${\bf{x}}$ and ${\bf{t}}$, and $\mathcal{L}(y; {\bf{x}}, {\bf{t}})$ is a covariate-dependent slowly varying function at infinity; i.e., $\mathcal{L}(ys; {\bf{x}}, {\bf{t}})/\mathcal{L}(y; {\bf{x}}, {\bf{t}})\to1$ as $y\to\infty$ for any $s>0$. We assume that the slowly varying function $\mathcal{L}(y; {\bf{x}}, {\bf{t}})$ converges to a constant at a reasonably high speed. Specifically, we assume
\begin{equation}
	\mathcal{L}(y; {\bf{x}}, {\bf{t}})=c_0({\bf{x}}, {\bf{t}})+c_1({\bf{x}}, {\bf{t}})y^{-\beta({\bf{x}}, {\bf{t}})}+o(y^{-\beta({\bf{x}}, {\bf{t}})}),\label{2.1.2}\tag{2.2}
\end{equation}
where $c_0({\bf{x}}, {\bf{t}})$, $c_1({\bf{x}}, {\bf{t}})$ and $\beta({\bf{x}}, {\bf{t}})$ are functions of ${\bf{x}}$ and ${\bf{t}}$ with $c_0({\bf{x}}, {\bf{t}})>0$ and $\beta({\bf{x}}, {\bf{t}})>0$, and $o(y^{-\beta({\bf{x}}, {\bf{t}})})$ is a higher-order term. This assumption is called the Hall class (Hall 1982). In our study, we adopt the varying coefficient model for the conditional extreme value index $\gamma({\bf{x}}, {\bf{t}})$ as
\begin{equation}
	\log\gamma({\bf{x}}, {\bf{t}})^{-1}=(1, {\bf{x}}^\top){\bm{\theta}}({\bf{t}})=\theta_0({\bf{t}})+\theta_1({\bf{t}})x_1+\cdots+\theta_p({\bf{t}})x_p,\label{2.1.3}\tag{2.3}
\end{equation}
where ${\bf{x}}=(x_1, x_2, \ldots, x_p)^\top\in\mathbb{R}^p$, ${\bm{\theta}}({\bf{t}})=(\theta_0({\bf{t}}), \theta_1({\bf{t}}), \ldots, \theta_p({\bf{t}}))^\top\in\mathbb{R}^{p+1}$, and $\theta_j({\bf{t}}),\ j=0, 1, \ldots, p$ are the unknown smooth functions of ${\bf{t}}$.

\subsection{Local constant maximum likelihood estimator}\label{subsec2.2}

Let $f(y\mid {\bf{x}}, {\bf{t}})$ be the conditional probability density function of $Y$ given $({\bf{X}}, {\bf{T}})=({\bf{x}}, {\bf{t}})$. If $\mathcal{L}(\cdot; {\bf{x}}, {\bf{t}})$ is differentiable, we obtain
\begin{equation*}
	f(y\mid {\bf{x}}, {\bf{t}})=\gamma({\bf{x}}, {\bf{t}})^{-1}y^{-1/\gamma({\bf{x}}, {\bf{t}})-1}\mathcal{L}(y; {\bf{x}}, {\bf{t}})-y^{-1/\gamma({\bf{x}}, {\bf{t}})}\frac{\partial\mathcal{L}}{\partial y}(y; {\bf{x}}, {\bf{t}}).
\end{equation*}
Because $\mathcal{L}(y; {\bf{x}}, {\bf{t}})\to c_0({\bf{x}}, {\bf{t}})$ and $\partial\mathcal{L}(y; {\bf{x}}, {\bf{t}})/\partial y\to0$ as $y\to\infty$ by using (\ref{2.1.2}), we have
\begin{equation*}
	f(y\mid {\bf{x}}, {\bf{t}})\approx\frac{c_0({\bf{x}}, {\bf{t}})}{\gamma({\bf{x}}, {\bf{t}})}y^{-1/\gamma({\bf{x}}, {\bf{t}})-1}
\end{equation*}
for sufficiently large $y>0$. Let $\{(Y_i, {\bf{X}}_i, {\bf{T}}_i),\ i=1, 2, \ldots, n\}$ be an independent and identically distributed random sample with the same distribution as $(Y, {\bf{X}}, {\bf{T}})$. We introduce a sufficiently high threshold $\omega_n>0$ such that $\omega_n\to\infty$ as $n\to\infty$ and use the responses that exceed it. Let $f(y\mid {\bf{x}}, {\bf{t}}, \omega_n)$ be the conditional probability density function of $Y$ given $({\bf{X}}, {\bf{T}})=({\bf{x}}, {\bf{t}})$ and $Y>\omega_n$. Then, we have
\begin{equation}
	f(y\mid {\bf{x}}, {\bf{t}}, \omega_n)\approx\gamma({\bf{x}}, {\bf{t}})^{-1}\left(\frac{y}{\omega_n}\right)^{-1/\gamma({\bf{x}}, {\bf{t}})}y^{-1}\label{2.2.4}\tag{2.4}
\end{equation}
for $y>\omega_n$. Thus, we can estimate the coefficient function vector ${\bm{\theta}}({\bf{t}})$ by using the following weighted maximum likelihood approach: Let
\begin{equation}
	L_n({\bm{\theta}})=\sum_{i=1}^n\left\{\exp((1, {\bf{X}}_i^\top){\bm{\theta}})\log\frac{Y_i}{\omega_n}-(1, {\bf{X}}_i^\top){\bm{\theta}}\right\}I(Y_i>\omega_n)K({\bf{H}}_n^{-1}({\bf{t}}-{\bf{T}}_i)),\label{2.2.5}\tag{2.5}
\end{equation}
where ${\bm{\theta}}\in\mathbb{R}^{p+1}$, $I(\cdot)$ is an indicator function, $K(\cdot)$ is a kernel function on $\mathbb{R}^q$, and ${\bf{H}}_n={\rm{diag}}(h_{n1}, \ldots, h_{nq})$ is a $q$-order diagonal matrix whose components are bandwidths $h_{nk},\ k=1, 2, \ldots, q$ such that $h_{nk}\to0$ as $n\to\infty$. We define the estimator of the coefficient function vector ${\bm{\theta}}({\bf{t}})$ as the minimizer of the objective function $L_n({\bm{\theta}})$. We denote this estimator by $\widehat{\bm{\theta}}({\bf{t}})=(\widehat{\theta}_0({\bf{t}}), \widehat{\theta}_1({\bf{t}}), \ldots, \widehat{\theta}_p({\bf{t}}))^\top\in\mathbb{R}^{p+1}$. Ma et al. (2019) provided the local linear maximum likelihood estimator. When $p=0$ and $q=0$, the covariate-independent estimator $\widehat{\theta}_0$ is explicitly obtained, and we have
\begin{equation*}
	\widehat{\gamma}=\exp(-\widehat{\theta}_0)=\frac{\sum_{i=1}^n(\log Y_i-\log\omega_n)I(Y_i>\omega_n)}{\sum_{i=1}^nI(Y_i>\omega_n)},
\end{equation*}
which is the Hill estimator proposed by Hill (1975) and is widely used in univariate extreme value theory.

Note that the varying coefficient model corresponds to linear and nonparametric models as special cases. When $q=0$, (\ref{2.1.3}) is simplified as
\begin{equation*}
	\log\gamma({\bf{x}})^{-1}=(1, {\bf{x}}^\top){\bm{\theta}}=\theta_0+\theta_1x_1+\cdots+\theta_px_p,
\end{equation*}
where ${\bm{\theta}}=(\theta_0, \theta_1, \ldots, \theta_p)^\top\in\mathbb{R}^{p+1}$, and $\theta_j,\ j=0, 1, \ldots, p$ are the regression coefficients. Wang and Tsai (2009) studied this tail index regression model. Whereas, when $p=0$, we obtain a nonparametric estimator of the positive extreme value index as
\begin{equation*}
	\widehat{\gamma}({\bf{t}})=\exp(-\widehat{\theta}_0({\bf{t}}))=\frac{\sum_{i=1}^n(\log Y_i-\log\omega_n)I(Y_i>\omega_n)K({\bf{H}}_n^{-1}({\bf{t}}-{\bf{T}}_i))}{\sum_{i=1}^nI(Y_i>\omega_n)K({\bf{H}}_n^{-1}({\bf{t}}-{\bf{T}}_i))},
\end{equation*}
which was studied by Goegebeur et al. (2014, 2015).

\subsection{Bandwidths and threshold selection}\label{subsec2.3}

The threshold $\omega_n$ and bandwidths $h_{nk},\ k=1, \ldots, q$ are tuning parameters that control the balance between the bias and variance of the estimator $\widehat{\bm{\theta}}({\bf{t}})$. A larger value of $h_{nk}$ or smaller value of $\omega_n$ leads to more bias, whereas a larger value of $\omega_n$ or smaller value of $h_{nk}$ leads to a larger variance. Therefore, these tuning parameters must be appropriately selected.

The threshold selection is needed to obtain the good approximation of (\ref{2.2.4}). To achieve this, the discrepancy measure, which was proposed by Wang and Tsai (2009), is suitable. Meanwhile, the choice of the bandwidths controls the smoothness of the estimator. Therefore, we use the cross-validation to select the bandwidths, similar to other studies using kernel smoothing (e.g., Ma et al. 2019). Thus, we combine the discrepancy measure and cross-validation as the whole tuning parameter selection method. The algorithm of the tuning parameter selection is as follows. In the first step, we select the bandwidths $h_{nk},\ k=1, \ldots, q$ by $D$-fold cross-validation as
\begin{align*}
	{\bf{H}}_{\rm{CV}}&:=\argmin_{\bf{H}}\sum_{d=1}^D\sum_{i=1}^{\lfloor n/D\rfloor}\left\{\exp((1, {\bf{X}}_i^{(d)}{}^\top)\widehat{\bm{\theta}}^{(d)}({\bf{T}}_i^{(d)}; \omega_0, {\bf{H}}))\log\frac{Y_i^{(d)}}{\omega_0}\right.\\
	&\hspace*{3cm}\Biggl.-(1, {\bf{X}}_i^{(d)}{}^\top)\widehat{\bm{\theta}}^{(d)}({\bf{T}}_i^{(d)}; \omega_0, {\bf{H}})\Biggr\}I(Y_i^{(d)}>\omega_0),
\end{align*}
which is based on (\ref{2.2.5}), where $\omega_0$ is a predetermined threshold, $\lfloor\cdot\rfloor$ is the floor function, and $\{(Y_i^{(d)}, {\bf{X}}_i^{(d)}, {\bf{T}}_i^{(d)}),\ i=1, 2, \ldots, \lfloor n/D\rfloor\}$ is the $d$th test dataset. $\widehat{\bm{\theta}}^{(d)}(\cdot; \omega_0, {\bf{H}})$ is the proposed estimator, with $\omega_n=\omega_0$ and ${\bf{H}}_n={\bf{H}}$, which is obtained from the $d$th training dataset. In the second step, we select the threshold $\omega_n$ using the discrepancy measure. We denote the order statistics of $\{\exp\{-\exp((1, {\bf{X}}_i^\top)\widehat{\bm{\theta}}({\bf{T}}_i))\log(Y_i/\omega_n)\}: Y_i>\omega_n,\ i=1, \ldots, n\}$ as $\widehat{U}_{1, n_0} \leq\widehat{U}_{2, n_0} \leq\ldots \leq\widehat{U}_{n_0, n_0}$, where $n_0=\sum_{i=1}^nI(Y_i>\omega_n)$ is the number of responses that exceed the threshold $\omega_n$. Because the conditional distribution of $\exp\{-\exp((1, {\bf{X}}^\top){\bm{\theta}}({\bf{T}}))\log(Y/\omega_n)\}$ given $Y>\omega_n$ is approximately standard uniform, we can regard $\{\widehat{U}_{l, n_0}\}_{l=1}^{n_0}$ as a sample from the standard uniform distribution. Therefore, we select the threshold $\omega_n$ as
\begin{equation*}
	\omega_{\rm{DM}}:=\argmin_\omega\frac{1}{n_0}\sum_{l=1}^{n_0}(\widehat{U}_{l, n_0}(\omega, {\bf{H}}_{\rm{CV}})-\widehat{F}(l/n_0; \omega, {\bf{H}}_{\rm{CV}}))^2,
\end{equation*}
where $\{\widehat{U}_{l, n_0}(\omega, {\bf{H}}_{\rm{CV}})\}_{l=1}^{n_0}$ are $\{\widehat{U}_{l, n_0}\}_{l=1}^{n_0}$ with $\omega_n=\omega$ and ${\bf{H}}_n={\bf{H}}_{\rm{CV}}$, and $\widehat{F}(\cdot; \omega, {\bf{H}}_{\rm{CV}})$ is a empirical distribution function of $\{\widehat{U}_{l, n_0}(\omega, {\bf{H}}_{\rm{CV}})\}_{l=1}^{n_0}$.

\section{Asymptotic properties}\label{sec3}

\subsection{Conditions}\label{subsec3.1}

In this section, we investigate the asymptotic properties of our proposed estimator. The following technical conditions are required: We define $n_0({\bf{t}})=n\det({\bf{H}}_n)f_{\bf{T}}({\bf{t}})P(Y>\omega_n\mid {\bf{T}}={\bf{t}})$ and $n({\bf{t}})=n\det({\bf{H}}_n)f_{\bf{T}}({\bf{t}})$, where $f_{\bf{T}}({\bf{t}})$ is the marginal probability density function of ${\bf{T}}$. We also define
\begin{equation*}
	{\bm{\Sigma}}_n({\bf{t}})=E[(1, {\bf{X}}^\top)^\top(1, {\bf{X}}^\top)\mid {\bf{T}}={\bf{t}}, Y>\omega_n]
\end{equation*}
and
\begin{equation*}
	\widetilde{\bm{\Sigma}}_n({\bf{t}})=\frac{1}{n_0({\bf{t}})}\sum_{i=1}^n(1, {\bf{X}}_i^\top)^\top(1, {\bf{X}}_i^\top)I(Y_i>\omega_n)K({\bf{H}}_n^{-1}({\bf{t}}-{\bf{T}}_i)).
\end{equation*}

\begin{itemize}
	\item[(C.1)]The kernel function $K(\cdot)$ is an absolutely continuous function that has compact support and satisfies the conditions
	\begin{equation*}
		\int K({\bf{u}}){\rm{d}}{\bf{u}}=1,\ \int u_kK({\bf{u}}){\rm{d}}{\bf{u}}=0,\ \int u_k^2K({\bf{u}}){\rm{d}}{\bf{u}}<\infty,\ \int K({\bf{u}})^2{\rm{d}}{\bf{u}}<\infty
	\end{equation*}
	with ${\bf{u}}=(u_1, u_2, \ldots, u_q)^\top\in\mathbb{R}^q$.
	\item[(C.2)]The joint probability density function $f(y, {\bf{x}}, {\bf{t}})$ of $(Y, {\bf{X}}, {\bf{T}})$ and the coefficient function $\theta_j({\bf{t}})$ have continuous second-order derivative on ${\bf{t}}$.
	\item[(C.3)]Assume $n_0({\bf{t}})\to\infty$ and
	\begin{equation*}
		{\bm{\Sigma}}_n({\bf{t}})^{-1/2}\widetilde{\bm{\Sigma}}_n({\bf{t}}){\bm{\Sigma}}_n({\bf{t}})^{-1/2}\xrightarrow{P}{\bf{I}}_{p+1}
	\end{equation*}
	as $n\to\infty$ for all ${\bf{t}}\in\mathbb{R}^q$, where ${\bf{I}}_{p+1}$ is a $(p+1)$-order identity matrix and the symbol ``$\xrightarrow{P}$'' stands for convergence in probability.
	\item[(C.4)]The reminder term $o(y^{-\beta({\bf{x}}, {\bf{t}})})$ defined in (\ref{2.1.2}) satisfies
	\begin{equation*}
		\sup_{{\bf{x}}\in\mathbb{R}^p,\ {\bf{t}}\in\mathbb{R}^q}(y^{\beta({\bf{x}}, {\bf{t}})}o(y^{-\beta({\bf{x}}, {\bf{t}})}))\to0
	\end{equation*}
	as $y\to\infty$.
	\item[(C.5)]For all ${\bf{t}}\in\mathbb{R}^q$, there exists a nonzero vector ${\bf{b}}({\bf{t}})\in\mathbb{R}^{p+1}$ such that
	\begin{align*}
		&\frac{n({\bf{t}})}{\sqrt{n_0({\bf{t}})}}{\bf{\Sigma}}_n({\bf{t}})^{-1/2}\\
		&E\left[(1, {\bf{X}}^\top)^\top c_1({\bf{X}}, {\bf{t}})\frac{\gamma({\bf{X}}, {\bf{t}})\beta({\bf{X}}, {\bf{t}})}{1+\gamma({\bf{X}}, {\bf{t}})\beta({\bf{X}}, {\bf{t}})}\omega_n^{-1/\gamma({\bf{X}}, {\bf{t}})-\beta({\bf{X}}, {\bf{t}})}\mid {\bf{T}}={\bf{t}}\right]\to{\bf{b}}({\bf{t}})
	\end{align*}
	as $n\to\infty$.
\end{itemize}
The condition (C.1) is typically used for kernel estimation. The conditions (C.3)-(C.5) correspond to the conditions (C.1)-(C.3) of Wang and Tsai (2009). The condition (C.3) requires that a certain weak law of large numbers holds. The condition (C.4) regularizes the extreme behavior of the slowly varying function $\mathcal{L}(y; {\bf{x}}, {\bf{t}})$. The condition (C.5) specifies the optimal convergence rates of threshold $\omega_n$ and bandwidths $h_{nk},\ k=1, \ldots, q$.

\subsection{Asymptotic properties}\label{subsec3.2}

We define
\begin{equation*}
	\dot{\bm{L}}_n({\bm{\theta}})=\sum_{i=1}^n(1, {\bf{X}}_i^\top)^\top\left\{\exp((1, {\bf{X}}_i^\top){\bm{\theta}})\log\frac{Y_i}{\omega_n}-1\right\}I(Y_i>\omega_n)K({\bf{H}}_n^{-1}({\bf{t}}-{\bf{T}}_i))
\end{equation*}
and
\begin{equation*}
	\ddot{\bm{L}}_n({\bm{\theta}})=\sum_{i=1}^n(1, {\bf{X}}_i^\top)^\top(1, {\bf{X}}_i^\top)\exp((1, {\bf{X}}_i^\top){\bm{\theta}})\log\frac{Y_i}{\omega_n}I(Y_i>\omega_n)K({\bf{H}}_n^{-1}({\bf{t}}-{\bf{T}}_i)).
\end{equation*}
The above $\dot{\bm{L}}_n({\bm{\theta}})$ and $\ddot{\bm{L}}_n({\bm{\theta}})$ are the gradient vector and Hessian matrix of the objective function $L_n({\bm{\theta}})$, respectively. The proposed estimator $\widehat{\bm{\theta}}({\bf{t}})$ is defined as the minimizer of $L_n({\bm{\theta}})$ and satisfies $\dot{\bm{L}}_n(\widehat{\bm{\theta}}({\bf{t}}))={\bm{0}}$. Therefore, similar to common approaches for establishing the asymptotic normality of the maximum likelihood estimator, we need to investigate the asymptotic properties of $\dot{\bm{L}}_n({\bm{\theta}})$ and $\ddot{\bm{L}}_n({\bm{\theta}})$. Let $\nu=\int K({\bf{u}})^2{\rm{d}}{\bf{u}}$, ${\bm{\kappa}}=\int {\bf{u}}{\bf{u}}^\top K({\bf{u}}){\rm{d}}{\bf{u}}$,
\begin{equation*}
	{\bm{\Delta}}_{k}({\bf{t}})=\left(\frac{\partial \theta_0}{\partial t_k}({\bf{t}}), \frac{\partial \theta_1}{\partial t_k}({\bf{t}}), \ldots, \frac{\partial \theta_p}{\partial t_k}({\bf{t}})\right)^\top
\end{equation*}
and 
\begin{equation*}
	{\bm{\Delta}}_{k_1k_2}({\bf{t}})=\left(\frac{\partial^2 \theta_0}{\partial t_{k_1}\partial t_{k_2}}({\bf{t}}), \frac{\partial^2 \theta_1}{\partial t_{k_1}\partial t_{k_2}}({\bf{t}}), \ldots,\frac{\partial^2 \theta_p}{\partial t_{k_1}\partial t_{k_2}}({\bf{t}})\right)^\top,
\end{equation*}
where ${\bf{t}}=(t_1, t_2, \ldots, t_q)^\top\in\mathbb{R}^q$ and $k_1, k_2\in\{1, 2, \ldots, q\}$.

\begin{theorem}\label{thm1}
	Let us suppose that conditions (C.1)-(C.5) are satisfied; then, as $n\to\infty$,
	\begin{equation*}
		[n_0({\bf{t}}){\bm{\Sigma}}_n({\bf{t}})]^{-1/2}\left\{\dot{\bm{L}}_n({\bm{\theta}}({\bf{t}}))+n_0({\bf{t}}){\bm{\Sigma}}_n({\bf{t}})\sum_{l=1}^2{\bf{\Lambda}}_n^{(l)}({\bf{t}})\right\}\xrightarrow{D}{\rm{N}}(-{\bf{b}}({\bf{t}}), \nu{\bf{I}}_{p+1}),
	\end{equation*}
	where the symbol ``$\xrightarrow{D}$'' denotes convergence in the distribution,
	\begin{equation*}
		{\bf{\Lambda}}_n^{(l)}({\bf{t}})=\frac{1}{2}{\bm{\Sigma}}_n({\bf{t}})^{-1}E[(1, {\bf{X}}^\top)^\top{\rm{tr}}({\bm{\lambda}}^{(l)}({\bf{X}}, {\bf{t}}){\bf{H}}_n{\bm{\kappa}}{\bf{H}}_n)\mid {\bf{T}}={\bf{t}}, Y>\omega_n]
	\end{equation*}
	and 
	\begin{equation*}
		{\bm{\lambda}}^{(l)}({\bf{X}}, {\bf{t}})=\begin{cases}
			({\bm{\Delta}}_1({\bf{t}}), \ldots, {\bm{\Delta}}_q({\bf{t}}))^\top(1, {\bf{X}}^\top)^\top(1, {\bf{X}}^\top)({\bm{\Delta}}_1({\bf{t}}), \ldots, {\bm{\Delta}}_q({\bf{t}})), & l=1,\\
			[(1,{\bf{X}}^\top){\bm{\Delta}}_{k_1k_2}({\bf{t}})]_{q\times q}, & l=2.
		\end{cases}
	\end{equation*}
\end{theorem}

\begin{theorem}\label{thm2}
	Let us suppose that conditions (C.1)-(C.5) are satisfied; then, as $n\to\infty$,
	\begin{equation*}
		n_0({\bf{t}})^{-1}{\bm{\Sigma}}_n({\bf{t}})^{-1/2}\ddot{\bm{L}}_n({\bm{\theta}}({\bf{t}})){\bm{\Sigma}}_n({\bf{t}})^{-1/2}\xrightarrow{P}{\bf{I}}_{p+1}.
	\end{equation*}
\end{theorem}

From Theorems \ref{thm1} and \ref{thm2}, we obtain the following asymptotic normality of our proposed estimator $\widehat{\bm{\theta}}({\bf{t}})$:

\begin{theorem}\label{thm3}
	Let us suppose that conditions (C.1)-(C.5) are satisfied; then, as $n\to\infty$,
	\begin{equation*}
		[n_0({\bf{t}}){\bm{\Sigma}}_n({\bf{t}})]^{1/2}\left\{(\widehat{\bm{\theta}}({\bf{t}})-{\bm{\theta}}({\bf{t}}))-\sum_{l=1}^2{\bf{\Lambda}}_n^{(l)}({\bf{t}})\right\}\xrightarrow{D}{\rm{N}}({\bf{b}}({\bf{t}}), \nu{\bf{I}}_{p+1}).
	\end{equation*}
\end{theorem}

This result implies that $\widehat{\bm{\theta}}({\bf{t}})$ is the consistent estimator of ${\bm{\theta}}({\bf{t}})$. The convergence rate of $\widehat{\bm{\theta}}({\bf{t}})$ to ${\bm{\theta}}({\bf{t}})$ is on the same order as $[n_0({\bf{t}}){\bm{\Sigma}}_n({\bf{t}})]^{-1/2}$. The $n_0({\bf{t}})=n\det({\bf{H}}_n)f_{\bf{T}}({\bf{t}})P(Y>\omega_n\mid{\bf{T}}={\bf{t}})$ is proportional to the number of top-order statistics of the responses used for estimation at ${\bf{t}}$. The ${\bm{\Sigma}}_n({\bf{t}})$ is defined in Section \ref{subsec3.1}. The asymptotic bias is caused by two factors. The bias ${\bf{b}}({\bf{t}})$ is caused by the approximation of the tail of the conditional distribution of $Y$ by the Pareto distribution in (\ref{2.2.4}), which is related to the convergence rate of the slowly varying function $\mathcal{L}(\cdot; {\bf{x}}, {\bf{t}})$ to the constant $c_0({\bf{x}}, {\bf{t}})$. From the definition of ${\bf{b}}({\bf{t}})$ given in (C.5), we can see that the proposed estimator is more biased for larger $\gamma({\bf{x}}, {\bf{t}})$. In other words, the heavier the tail of the data, the more biased the estimator. Meanwhile, if $\beta({\bf{x}}, {\bf{t}})$ is small, the large bias of the estimator is occurred. Thus, the bias of our estimator is particularly sensitive to $\gamma({\bf{x}}, {\bf{t}})$ and $\beta({\bf{x}}, {\bf{t}})$. These parameters are related to the second order condition in extreme value theory (see, Theorem 3.2.5 of de Haan and Ferreira 2006, Theorems 2 and 3 of Wang and Tsai 2009, and Theorem 2 of Li et al. 2022, to name a few). In contrast, the biases ${\bf{\Lambda}}_n^{(1)}({\bf{t}})$ and ${\bf{\Lambda}}_n^{(2)}({\bf{t}})$ are caused by kernel smoothing.

Our asymptotic normality is comparable to the asymptotic normality of the local linear maximum likelihood estimator of the coefficient function vector proposed by Ma et al. (2019). The difference between the two estimators is the asymptotic bias. In the asymptotic normality in Ma et al. (2019), it is assumed that the bias caused by the approximation (\ref{2.2.4}) is negligible, so the bias ${\bf{b}}({\bf{t}})$ does not appear in their asymptotic normality. The essential difference is the bias caused by kernel smoothing. In the case of Ma et al. (2019), the bias caused by kernel smoothing is ${\bf{\Lambda}}_n^{(2)}({\bf{t}})$. However, it has the same convergence rate as the bias ${\bf{\Lambda}}_n^{(1)}({\bf{t}})+{\bf{\Lambda}}_n^{(2)}({\bf{t}})$. The asymptotic variances of the two estimators are the same.

\section{Testing for structure of the coefficient function}\label{sec4}

\subsection{Testing method}\label{subsec4.1}

In varying coefficient models, we often hope to test whether each coefficient function $\theta_j(\cdot)$ is constant or zero. If some $\theta_j({\bf{t}})$ does not vary across ${\bf{t}}$, this motivates us to consider models that are simpler than the varying coefficient model. Generally, the hypothesis test can be represented as
\begin{equation}
	{\rm{H}}_0: \theta_j(\cdot)\equiv\eta(\cdot)\quad{\rm{vs.}}\quad{\rm{H}}_1: \theta_j(\cdot)\not\equiv\eta(\cdot)\label{4.1.1}\tag{4.1}
\end{equation}
for a given known function $\eta(\cdot)$, where ${\rm{H}}_0$ is the null hypothesis and ${\rm{H}}_1$ is the alternative hypothesis.

Without a loss of generality, we assume that the explanatory variable vector ${\bf{T}}\in\mathbb{R}^q$ is distributed on $[0, 1]^q\subset\mathbb{R}^q$. Then, we apply Lemma 1 of Fan and Zhang (2000) to
\begin{equation*}
	\psi({\bf{t}})=\frac{1}{\sqrt{n({\bf{t}})\sigma_{nj}({\bf{t}})}}\left.\frac{\partial L_n({\bm{\theta}})}{\partial\theta_j}\mid_{{\bm{\theta}}=(\widehat{\theta}_0({\bf{t}}), \ldots, \widehat{\theta}_{j-1}({\bf{t}}), \theta_j({\bf{t}}), \widehat{\theta}_{j+1}({\bf{t}}), \ldots, \widehat{\theta}_p({\bf{t}}))^\top},\right.
\end{equation*}
where $\sigma_{nj}({\bf{t}})=E[X_j^2I(Y>\omega_n)\mid {\bf{T}}={\bf{t}}],\ j=0, 1, \ldots, p$ and $X_0\equiv1$. The following conditions are required:

\begin{itemize}
	\item[(C.6)]For all large $n\in\mathbb{N}$, the function $\sigma_{nj}({\bf{t}})$ is bounded away from zero for all ${\bf{t}}\in[0, 1]^q$ and has a bounded partial derivative.
	\item[(C.7)]$\lim_{n\to\infty}\sup_{\bf{t}}E[\lvert X_j\rvert^sI(Y>\omega_n)\mid {\bf{T}}={\bf{t}}]<\infty$ for some $s>2$.
\end{itemize}

\begin{theorem}\label{thm4}
	Under the conditions (C.1)-(C.7), if $h_n:=h_{nk}=n^{-b},\ k=1, \ldots, q$, for some $0<b<1-2/s$, we have
	\begin{equation*}
		P\left(\sqrt{-2q\log h_n}\left(\frac{1}{\sqrt{\nu}}\sup_{{\bf{t}}\in[0, 1]^q}\lvert\psi({\bf{t}})-E[\psi({\bf{t}})]\rvert-d_n\right)<s\right)\xrightarrow{D}\exp(-2\exp(-s))
	\end{equation*}
	as $n\to\infty$, where 
	\begin{align*}
		d_n&=\sqrt{-2q\log h_n}\\
		&\quad+\frac{1}{\sqrt{-2q\log h_n}}\left[\frac{1}{2}(q-1)\log\log(h_n^{-1})+\log\left\{\left(\frac{2q}{\pi}\right)^{q/2}\sqrt{\frac{{\rm{det}}({{\bm{\Xi}}})}{4q\pi}}\right\}\right]
	\end{align*}
	(Rosenblatt 1976) and 
	\begin{equation*}
		{\bm{\Xi}}=\frac{1}{2\nu}\left[\int\frac{\partial^2 K({\bf{u}})}{\partial u_{k_1}\partial u_{k_2}}{\rm{d}}{\bf{u}}\right]_{q\times q}.
	\end{equation*}
\end{theorem}

From Theorem \ref{thm3}, we now have $\widehat{\bm{\theta}}({\bf{t}})\to^P{\bm{\theta}}({\bf{t}})$ as $n\to\infty$. By the first-order Taylor expansion around $\widehat{\theta}_j({\bf{t}})=\theta_j({\bf{t}})$, we obtain
\begin{align*}
	&\frac{\partial L_n({\bm{\theta}})}{\partial\theta_j}\mid_{{\bm{\theta}}=(\widehat{\theta}_0({\bf{t}}), \widehat{\theta}_{1}({\bf{t}}),  \ldots, \widehat{\theta}_p({\bf{t}}))^\top}\\
	&=\frac{\partial L_n({\bm{\theta}})}{\partial\theta_j}\mid_{{\bm{\theta}}=(\widehat{\theta}_0({\bf{t}}), \ldots, \widehat{\theta}_{j-1}({\bf{t}}), \theta_j({\bf{t}}), \widehat{\theta}_{j+1}({\bf{t}}), \ldots, \widehat{\theta}_p({\bf{t}}))^\top}\\
	&\quad+\frac{\partial^2 L_n({\bm{\theta}})}{\partial\theta_j^2}\mid_{{\bm{\theta}}=(\widehat{\theta}_0({\bf{t}}), \ldots, \widehat{\theta}_{j-1}({\bf{t}}), \theta_j({\bf{t}}), \widehat{\theta}_{j+1}({\bf{t}}), \ldots, \widehat{\theta}_p({\bf{t}}))^\top}(\widehat{\theta}_j({\bf{t}})-\theta_j({\bf{t}}))\{1+o_P(1)\}.
\end{align*}
The left-hand side of the above equation is zero because $\widehat{\bm{\theta}}({\bf{t}})=(\widehat{\theta}_0({\bf{t}}), \widehat{\theta}_{1}({\bf{t}}),  \ldots, \widehat{\theta}_p({\bf{t}}))^\top$ is the minimizer of $L_n({\bm{\theta}})$. From Theorems \ref{thm2} and \ref{thm3}, we also have
\begin{equation*}
	\partial^2L_n({\bm{\theta}})/\partial\theta_j^2\mid_{{\bm{\theta}}=(\widehat{\theta}_0({\bf{t}}), \ldots, \widehat{\theta}_{j-1}({\bf{t}}), \theta_j({\bf{t}}), \widehat{\theta}_{j+1}({\bf{t}}), \ldots, \widehat{\theta}_p({\bf{t}}))^\top}\to^Pn({\bf{t}})\sigma_{nj}({\bf{t}})
\end{equation*}
as $n\to\infty$. Consequently, we have
\begin{align*}
	&\frac{\partial L_n({\bm{\theta}})}{\partial \theta_j}\mid_{{\bm{\theta}}=(\widehat{\theta}_0({\bf{t}}), \ldots, \widehat{\theta}_{j-1}({\bf{t}}), \theta_j({\bf{t}}), \widehat{\theta}_{j+1}({\bf{t}}), \ldots, \widehat{\theta}_p({\bf{t}}))^\top}\\
	&\hspace*{1cm}=-n({\bf{t}})\sigma_{nj}({\bf{t}})(\widehat{\theta}_j({\bf{t}})-\theta_j({\bf{t}}))\{1+o_P(1)\}.
\end{align*}
This implies that $\psi({\bf{t}})$ in Theorem \ref{thm4} can be approximated as
\begin{equation*}
	\psi({\bf{t}})\approx-\sqrt{n({\bf{t}})\sigma_{nj}({\bf{t}})}(\widehat{\theta}_j({\bf{t}})-\theta_j({\bf{t}})).
\end{equation*}
From this result, $E[\psi({\bf{t}})]$ is asymptotically equivalent to the $j$th component of $-{\bf{b}}({\bf{t}})-[n_0({\bf{t}}){\bm{\Sigma}}_n({\bf{t}})]^{1/2}\sum_{l=1}^2{\bf{\Lambda}}_n^{(l)}({\bf{t}})$. This bias involves many unknown parameters. In particular, $\beta({\bf{x}}, {\bf{t}})$ included in ${\bf{b}}({\bf{t}})$ corresponds to the so-called second order parameter (see, Gomes et al. 2002). However, the estimation method of the second order parameter has not yet been developed in the context of the tail index regression. Thus, at the present stage, checking that (C.5) is satisfied and estimating $E[\psi({\bf{t}})]$ are challenging and are posited as future work. Therefore, in this paper, we assume that $E[\psi({\bf{t}})]$ is zero, similar to Wang and Tsai (2009). Then, Theorem \ref{thm4} can be used to test if (\ref{4.1.1}). Under the null hypothesis ${\rm{H}}_0: \theta_j({\bf{t}})\equiv\eta({\bf{t}})$, we use the test statistic
\begin{equation*}
	\widetilde{T}=\sqrt{-2q\log h_n}\left\{\frac{1}{\sqrt{\nu}}\max_{\bf{t}}\left(\widehat{[n({\bf{t}})\sigma_{nj}({\bf{t}})]}^{1/2}\lvert\widehat{\theta}_j({\bf{t}})-\eta({\bf{t}})\rvert\right)-d_n\right\},
\end{equation*}
where $\widehat{[n({\bf{t}})\sigma_{nj}({\bf{t}})]}$ is the kernel estimator of $n({\bf{t}})\sigma_{nj}({\bf{t}})$ based on (C.3). For a given significance level $\alpha$, we reject the null hypothesis ${\rm{H}}_0$ if $\widetilde{T}<-\log\{-0.5\log(\alpha/2)\}=e_{\alpha/2}$ or $\widetilde{T}>-\log\{-0.5\log(1-\alpha/2)\}=e_{1-\alpha/2}$.

As mentioned above, we are mainly interested in the following two hypothesis tests: One is
\begin{equation*}
	{\rm{H}}_{0\rm{Z}}: \theta_j(\cdot)\equiv0\quad{\rm{vs.}}\quad{\rm{H}}_{1\rm{Z}}: \theta_j(\cdot)\not\equiv0.
\end{equation*}
If the null hypothesis ${\rm{H}}_{0\rm{Z}}$ is not rejected, the corresponding $X_j$ may not be important for predicting the tail behavior of the distribution of the response $Y$. Thus, this can help judge the sparsity of a particular covariate. The other is
\begin{equation*}
	{\rm{H}}_{0{\rm{C}}}: \theta_j(\cdot)\equiv C_0\quad{\rm{vs.}}\quad{\rm{H}}_{1{\rm{C}}}: \theta_j(\cdot)\not\equiv C_0
\end{equation*}
for an unknown constant $C_0$ without prior knowledge. Under the null hypothesis ${\rm{H}}_{0{\rm{C}}}$, we estimate the unknown constant $C_0$ as the average of the estimates $\{\widehat{\theta}_j({\bf{t}}_l)\}_{l=1}^L$, where ${\bf{t}}_l,\ l=1, 2, \ldots, L$ are equally spaced points in $[0, 1]^q$. If the null hypothesis ${\rm{H}}_{0\rm{C}}$ is not rejected, it motivates us to adopt a simpler model that considers the coefficient function $\theta_j(\cdot)$ to be constant.

The simultaneous test from Theorem \ref{thm4} is more rigorous than the test statistic based on the residual sum of squares (see, Cai et al. 2000). Here, we consider the separate hypotheses for each coefficient. One might think that the single hypothesis test on all coefficients would be of interest. However, such an extension is really difficult because we have to consider the distribution of $\sup_{\bf{t}}\{\widehat{\theta}_0({\bf{t}}), \widehat{\theta}_1({\bf{t}}), \ldots, \widehat{\theta}_p({\bf{t}})\}$. In fact, such a method has not even been studied in the context of mean regression. Thus, the development of a simultaneous testing method into a single hypothesis test on all coefficient functions is posited as an important future work.

\subsection{Simulation}\label{subsec4.2}

We ran a simulation study to demonstrate the finite sample performance of the proposed estimator and test statistic. We present the results for the three model settings. In all settings, we simulated the responses $\{Y_i\}_{i=1}^n$ using the following conditional distribution function:
\begin{align*}
	1-F(y\mid {\bf{x}} ,{\bf{t}})&=\frac{(1+\delta)y^{-1/\gamma({\bf{x}}, {\bf{t}})}}{1+\delta y^{-1/\gamma({\bf{x}}, {\bf{t}})}}\\
	&=y^{-1/\gamma({\bf{x}}, {\bf{t}})}\{(1+\delta)-\delta(1+\delta)y^{-1/\gamma({\bf{x}}, {\bf{t}})}+o(y^{-1/\gamma({\bf{x}}, {\bf{t}})})\},
\end{align*}
where $\log\gamma({\bf{x}}, {\bf{t}})^{-1}=(1, {\bf{x}}^\top){\bm{\theta}}({\bf{t}})$. This conditional distribution function satisfies (\ref{2.1.2}) with $c_0({\bf{x}},{\bf{t}})=1+\delta$, $c_1({\bf{x}}, {\bf{t}})=-\delta(1+\delta)$, and $\beta({\bf{x}}, {\bf{t}})=\gamma({\bf{x}}, {\bf{t}})^{-1}$. If $\delta\neq0$, the above conditional distribution is not the Pareto distribution; therefore, we need to introduce the threshold $\omega_n$ appropriately. Otherwise, modeling bias occurs, resulting in less accuracy in the estimation. We simulated the predictors $\{(X_{i1}, X_{i2}, \ldots, X_{ip})\}_{i=1}^n$ based on the following procedure:
\begin{equation*}
	X_{ij}=\sqrt{12}(\Phi(Z_{ij})-1/2),
\end{equation*}
where $\{(Z_{i1}, Z_{i2}, \ldots, Z_{ip})\}_{i=1}^n$ is an independent sample from the multivariate normal distribution with $E[Z_{ij}]=0$ and ${\rm{cov}}[Z_{ij_1}, Z_{ij_2}]=0.5^{\lvert j_1-j_2\rvert}$, and $\Phi(\cdot)$ is the cumulative distribution function of a standard normal. Consequently, for $j=1, 2, \ldots, p$, $\{X_{ij}\}_{i=1}^n$ is uniformly distributed on $[-\sqrt{3}, \sqrt{3}]$ with unit variance. Meanwhile, we simulated the predictors $\{(T_{i1}, T_{i2}, \ldots, T_{iq})\}_{i=1}^n$ from a uniform distribution on $[-0.2, 1.2]^q\subset\mathbb{R}^q$ with ${\rm{cov}}[T_{ik_1}, T_{ik_2}]=0$.

To measure the goodness of the estimator $\widehat{\bm{\theta}}({\bf{t}})$, we calculated the following empirical mean square error based on $M=100$ simulations:
\begin{equation*}
	{\rm{MSE}}=\frac{1}{L}\frac{1}{M}\sum_{l=1}^L\sum_{m=1}^M(\widehat{\theta}_j^{(m)}({\bf{t}}_l)-\theta_j({\bf{t}}_l))^2,
\end{equation*}
where $\widehat{\theta}_j^{(m)}(\cdot)$ is the estimate of $\theta_j(\cdot)$ using the $m$th dataset and $\{{\bf{t}}_l\}_{l=1}^L$ are equally spaced points in $[0, 1]^q$. In addition, to evaluate the performance of the test statistic, we obtained the probability of error as follows. When the null hypothesis is true, the empirical probability of the Type I error is defined as
\begin{equation*}
	{\rm{E1}}=\#\{m: \widetilde{T}_m<e_{\alpha/2}\ {\rm{or}}\ \widetilde{T}_m>e_{1-\alpha/2},\ m=1, 2, \ldots, M\}/M,
\end{equation*}
where $\widetilde{T}_m$ is the test statistic $\widetilde{T}$ using the $m$th dataset. Meanwhile, when the null hypothesis is false, the empirical probability of the Type II error is given by 
\begin{equation*}
	{\rm{E2}}=\#\{m: e_{\alpha/2}<\widetilde{T}_m<e_{1-\alpha/2},\ m=1, 2 , \ldots, M\}/M.
\end{equation*}
Now, the null hypotheses of interest, ${\rm{H}}_{0{\rm{C}}}$ and ${\rm{H}}_{0{\rm{Z}}}$, are defined in Section \ref{subsec4.1}. Accordingly, if the null hypothesis ${\rm{H}}_{0{\rm{C}}}$ is true, i.e., the given coefficient function $\theta_j({\bf{t}})$ is constant, we provide E1 to examine the performance of the constancy test; if not, E2 is provided. Similarly, if the null hypothesis ${\rm{H}}_{0{\rm{Z}}}$ is true, i.e., $\theta_j({\bf{t}})\equiv0$, E1 is used to evaluate the accuracy for the sparsity test; if not, the result for ${\rm{H}}_{0{\rm{Z}}}$ is given as E2.

\begin{table}[t]
	\begin{center}
		\caption{\footnotesize{Results of estimation and hypothesis testing in the first model setting. E1 and E2 represent the empirical probabilities of Type I and Type II errors, respectively. Since $\theta_1(t)=1$ is a nonzero constant, the null hypothesis ${\rm{H}}_{0{\rm{C}}}$ is true, whereas the null hypothesis ${\rm{H}}_{0{\rm{Z}}}$ is false. Accordingly, we provide E1 for ${\rm{H}}_{0{\rm{C}}}$ and E2 for ${\rm{H}}_{0{\rm{Z}}}$. Meanwhile, since $\theta_2(t)=\cos(2t)$ is not constant, the null hypotheses ${\rm{H}}_{0{\rm{C}}}$ and ${\rm{H}}_{0{\rm{Z}}}$ are false and thus E2 is given for both tests. Similarly, $\theta_3(t)=0$ indicates that the null hypotheses ${\rm{H}}_{0{\rm{C}}}$ and ${\rm{H}}_{0{\rm{Z}}}$ are true. Thus, for $\theta_3(t)=0$, E1 is provided for both tests.\\}}
		\label{tab1}
		\begin{tabular}{@{}rr|rrr|rrr|rrr@{}}
			\hline
			& & \multicolumn{3}{c|}{$\theta_1(t)=1$} & \multicolumn{3}{c|}{$\theta_2(t)=\cos(2t)$} & \multicolumn{3}{c}{$\theta_3(t)=0$}\\
			\hline
			& & & \multicolumn{1}{c}{${\rm{H}}_{0{\rm{C}}}$} & \multicolumn{1}{c|}{${\rm{H}}_{0{\rm{Z}}}$} & & \multicolumn{1}{c}{${\rm{H}}_{0{\rm{C}}}$} & \multicolumn{1}{c|}{${\rm{H}}_{0{\rm{Z}}}$} & & \multicolumn{1}{c}{${\rm{H}}_{0{\rm{C}}}$} & \multicolumn{1}{c}{${\rm{H}}_{0{\rm{Z}}}$}\\
			\multicolumn{1}{c}{$\delta$} & \multicolumn{1}{c|}{$n$} & \multicolumn{1}{c}{MSE} & \multicolumn{1}{c}{E1} & \multicolumn{1}{c|}{E2} & \multicolumn{1}{c}{MSE} & \multicolumn{1}{c}{E2} & \multicolumn{1}{c|}{E2} & \multicolumn{1}{c}{MSE} & \multicolumn{1}{c}{E1} & \multicolumn{1}{c}{E1}\\
			\hline\hline
			0.1 & 200 & 0.1591 & 0.00 & 0.00 & 0.2392 & 0.42 & 0.03 & 0.1733 & 0.01 & 0.04\\
			& 500 & 0.1083 & 0.00 & 0.00 & 0.1535 & 0.07 & 0.00 & 0.1076 & 0.00 & 0.03\\
			& 1000 & 0.0960 & 0.00 & 0.00 & 0.1157 & 0.00 & 0.00 & 0.0878 & 0.02 & 0.03\\\hline
			0.25 & 200 & 0.1769 & 0.01 & 0.00 & 0.2499 & 0.45 & 0.07 & 0.1967 & 0.02 & 0.07\\
			& 500 & 0.1462 & 0.01 & 0.00 & 0.1670 & 0.07 & 0.01 & 0.1311 & 0.03 & 0.01\\
			& 1000 & 0.1258 & 0.01 & 0.00 & 0.1282 & 0.00 & 0.00 & 0.1162 & 0.02 & 0.06\\\hline
			0.5 & 200 & 0.2217 & 0.01 & 0.00 & 0.2650 & 0.55 & 0.09 & 0.1864 & 0.00 & 0.01\\
			& 500 & 0.1802 & 0.00 & 0.00 & 0.1929 & 0.09 & 0.00 & 0.1341 & 0.00 & 0.00\\
			& 1000 & 0.1591 & 0.04 & 0.00 & 0.1428 & 0.00 & 0.00 & 0.1195 & 0.02 & 0.02\\
			\hline
		\end{tabular}
	\end{center}
\end{table}

In the first model setting, we set $p=3$ and $q=1$ and defined the coefficient functions $\theta_j(\cdot),\ j=1, 2, 3$ as
\begin{equation*}
	\theta_j(t)=\begin{cases}
		1, & j=1,\\
		\cos(2t), & j=2,\\
		0, & j=3,
	\end{cases}
\end{equation*}
where the intercept term $\theta_0(t)$ was not considered. We employed the Epanechnikov kernel in the proposed estimator. In the estimation process, we selected the threshold $\omega_n$ and bandwidth $h_n$ using the procedure described in Section \ref{subsec2.3}. We set the pre-determined sample fraction to $n_0/n=0.2$ in $D=20$-fold cross-validation, where $n_0=\sum_{i=1}^nI(Y_i>\omega_n)$. Table \ref{tab1} shows the calculated MSEs and empirical probabilities of error for each coefficient function $\theta_j(\cdot)$ when $\delta=0.1, 0.25, 0.5$ and $n=200, 500, 1000$. For each coefficient function $\theta_j(\cdot)$, the calculated MSEs improved as $n$ increased. This result is desirable and suggests the consistency of the proposed estimator. Note that when testing the null hypothesis ${\rm{H}}_{0{\rm{C}}}$, we must estimate the unknown constant $C_0$. Since the maximum deviation between $\widehat{\theta}_j(t)$ and the estimated $C_0$ tends to be smaller than the maximum deviation between $\widehat{\theta}_j(t)$ and the true value $C_0$, the empirical probabilities of the Type I error were smaller for the null hypothesis ${\rm{H}}_{0{\rm{C}}}$ than for the null hypothesis ${\rm{H}}_{0{\rm{Z}}}$. In all settings, the empirical probability of the Type II error improved as $n$ increased.

\begin{figure}[t!]
	\centering
	\includegraphics[keepaspectratio, width=17.5cm,angle=270]{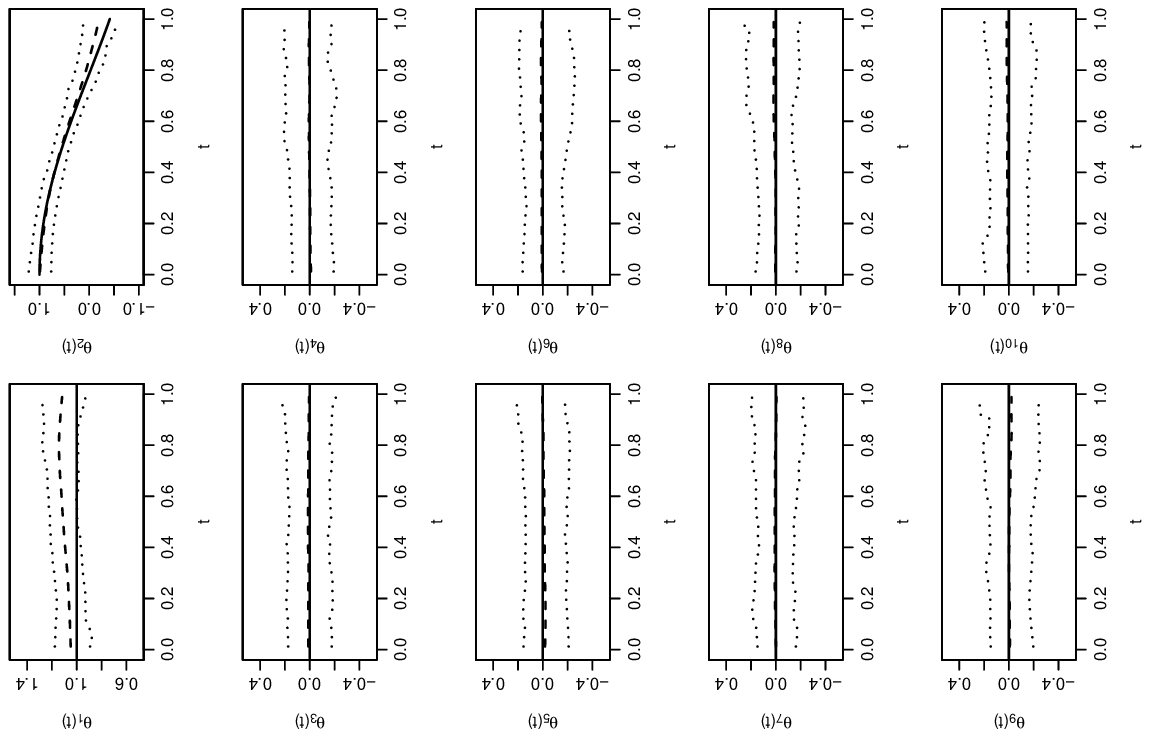}
	\caption{The estimated coefficient functions in the second model setting with $\delta=0.5$ and $n=1000$: the true value (solid line), average estimates (dashed line) and empirical 95\% confidence interval calculated based on estimates (dotted lines).}
	\label{fig1}
\end{figure}

\begin{table}[t!]
	\begin{center}
		\caption{\footnotesize{Results of estimation and hypothesis testing in the second model setting. E1 and E2 represent the empirical probabilities of Type I and Type II errors, respectively. Since $\theta_1(t)=1$ is a nonzero constant, the null hypothesis ${\rm{H}}_{0{\rm{C}}}$ is true, whereas the null hypothesis ${\rm{H}}_{0{\rm{Z}}}$ is false. Accordingly, we provide E1 for ${\rm{H}}_{0{\rm{C}}}$ and E2 for ${\rm{H}}_{0{\rm{Z}}}$. Meanwhile, since $\theta_2(t)=\cos(2t)$ is not constant, the null hypotheses ${\rm{H}}_{0{\rm{C}}}$ and ${\rm{H}}_{0{\rm{Z}}}$ are false and thus E2 is given for both tests. Similarly, $\theta_j(t)=0,\ j=3, 4, \ldots, 10$ indicate that the null hypotheses ${\rm{H}}_{0{\rm{C}}}$ and ${\rm{H}}_{0{\rm{Z}}}$ are true. Thus, for each $\theta_j(t)=0,\ j=3, 4, \ldots, 10$, E1 is provided for both tests.\\}}\label{tab2}
		\scalebox{0.73}{
			\begin{tabular}{@{}rr|rrr|rrr|rrr@{}}
				\hline
				& & \multicolumn{3}{c|}{$\theta_1(t)=1$} & \multicolumn{3}{c|}{$\theta_2(t)=\cos(2t)$} & \multicolumn{3}{c}{$\theta_3(t)=0$}\\
				\hline
				& & & \multicolumn{1}{c}{${\rm{H}}_{0{\rm{C}}}$} & \multicolumn{1}{c|}{${\rm{H}}_{0{\rm{Z}}}$} & & \multicolumn{1}{c}{${\rm{H}}_{0{\rm{C}}}$} & \multicolumn{1}{c|}{${\rm{H}}_{0{\rm{Z}}}$} & & \multicolumn{1}{c}{${\rm{H}}_{0{\rm{C}}}$} & \multicolumn{1}{c}{${\rm{H}}_{0{\rm{Z}}}$}\\
				\multicolumn{1}{c}{$\delta$} & \multicolumn{1}{c|}{$n$} & \multicolumn{1}{c}{MSE} & \multicolumn{1}{c}{E1} & \multicolumn{1}{c|}{E2} & \multicolumn{1}{c}{MSE} & \multicolumn{1}{c}{E2} & \multicolumn{1}{c|}{E2} & \multicolumn{1}{c}{MSE} & \multicolumn{1}{c}{E1} & \multicolumn{1}{c}{E1}\\
				\hline\hline
				0.1 & 200 & 0.1914 & 0.02 & 0.00 & 0.2640 & 0.64 & 0.07 & 0.1775 & 0.02 & 0.03\\
				& 500 & 0.1053 & 0.00 & 0.00 & 0.1945 & 0.12 & 0.00 & 0.1051 & 0.00 & 0.01\\
				& 1000 & 0.0802 & 0.00 & 0.00 & 0.1397 & 0.01 & 0.00 & 0.0782 & 0.00 & 0.03\\\hline
				0.25 & 200 & 0.2099 & 0.03 & 0.01 & 0.2717 & 0.56 & 0.10 & 0.2354 & 0.03 & 0.06\\
				& 500 & 0.1356 & 0.01 & 0.00 & 0.2003 & 0.29 & 0.00 & 0.1157 & 0.00 & 0.01\\
				& 1000 & 0.1125 & 0.00 & 0.00 & 0.1415 & 0.02 & 0.00 & 0.0915 & 0.00 & 0.02\\\hline
				0.5 & 200 & 0.2393 & 0.00 & 0.00 & 0.2818 & 0.66 & 0.15 & 0.2030 & 0.01 & 0.03\\
				& 500 & 0.1666 & 0.00 & 0.00 & 0.2203 & 0.30 & 0.01 & 0.1345 & 0.00 & 0.05\\
				& 1000 & 0.1333 & 0.01 & 0.00 & 0.1551 & 0.04 & 0.00 & 0.1071 & 0.01 & 0.04\\
				\hline\hline
				& & \multicolumn{3}{c|}{$\theta_4(t)=0$} & \multicolumn{3}{c|}{$\theta_5(t)=0$} & \multicolumn{3}{c}{$\theta_6(t)=0$}\\
				\hline
				& & & \multicolumn{1}{c}{${\rm{H}}_{0{\rm{C}}}$} & \multicolumn{1}{c|}{${\rm{H}}_{0{\rm{Z}}}$} & & \multicolumn{1}{c}{${\rm{H}}_{0{\rm{C}}}$} & \multicolumn{1}{c|}{${\rm{H}}_{0{\rm{Z}}}$} & & \multicolumn{1}{c}{${\rm{H}}_{0{\rm{C}}}$} & \multicolumn{1}{c}{${\rm{H}}_{0{\rm{Z}}}$}\\
				\multicolumn{1}{c}{$\delta$} & \multicolumn{1}{c|}{$n$} & \multicolumn{1}{c}{MSE} & \multicolumn{1}{c}{E1} & \multicolumn{1}{c|}{E1} & \multicolumn{1}{c}{MSE} & \multicolumn{1}{c}{E1} & \multicolumn{1}{c|}{E1} & \multicolumn{1}{c}{MSE} & \multicolumn{1}{c}{E1} & \multicolumn{1}{c}{E1}\\
				\hline\hline
				0.1 & 200 & 0.1911 & 0.02 & 0.03 & 0.1919 & 0.02 & 0.02 & 0.1786 & 0.03 & 0.04\\
				& 500 & 0.1061 & 0.00 & 0.01 & 0.1059 & 0.00 & 0.00 & 0.1058 & 0.00 & 0.03\\
				& 1000 & 0.0753 & 0.00 & 0.00 & 0.0783 & 0.00 & 0.03 & 0.0769 & 0.01 & 0.04\\\hline
				0.25 & 200 & 0.2188 & 0.03 & 0.05 & 0.1821 & 0.04 & 0.03 & 0.2110 & 0.04 & 0.04\\
				& 500 & 0.1263 & 0.00 & 0.05 & 0.1265 & 0.00 & 0.03 & 0.1342 & 0.00 & 0.03\\
				& 1000 & 0.1041 & 0.00 & 0.07 & 0.1020 & 0.00 & 0.04 & 0.0867 & 0.00 & 0.02\\\hline
				0.5 & 200 & 0.2074 & 0.02 & 0.06 & 0.2106 & 0.01 & 0.02 & 0.1906 & 0.00 & 0.03\\
				& 500 & 0.1409 & 0.00 & 0.05 & 0.1326 & 0.00 & 0.02 & 0.1423 & 0.01 & 0.03\\
				& 1000 & 0.1098 & 0.01 & 0.04 & 0.1085 & 0.00 & 0.01 & 0.1094 & 0.00 & 0.05\\
				\hline\hline
				& & \multicolumn{3}{c|}{$\theta_7(t)=0$} & \multicolumn{3}{c|}{$\theta_8(t)=0$} & \multicolumn{3}{c}{$\theta_9(t)=0$}\\
				\hline
				& & & \multicolumn{1}{c}{${\rm{H}}_{0{\rm{C}}}$} & \multicolumn{1}{c|}{${\rm{H}}_{0{\rm{Z}}}$} & & \multicolumn{1}{c}{${\rm{H}}_{0{\rm{C}}}$} & \multicolumn{1}{c|}{${\rm{H}}_{0{\rm{Z}}}$} & & \multicolumn{1}{c}{${\rm{H}}_{0{\rm{C}}}$} & \multicolumn{1}{c}{${\rm{H}}_{0{\rm{Z}}}$}\\
				\multicolumn{1}{c}{$\delta$} & \multicolumn{1}{c|}{$n$} & \multicolumn{1}{c}{MSE} & \multicolumn{1}{c}{E1} & \multicolumn{1}{c|}{E1} & \multicolumn{1}{c}{MSE} & \multicolumn{1}{c}{E1} & \multicolumn{1}{c|}{E1} & \multicolumn{1}{c}{MSE} & \multicolumn{1}{c}{E1} & \multicolumn{1}{c}{E1}\\
				\hline\hline
				0.1 & 200 & 0.1934 & 0.02 & 0.02 & 0.1868 & 0.03 & 0.03 & 0.1991 & 0.02 & 0.06\\
				& 500 & 0.0994 & 0.00 & 0.03 & 0.1103 & 0.01 & 0.01 & 0.1095 & 0.00 & 0.02\\
				& 1000 & 0.0853 & 0.00 & 0.06 & 0.0845 & 0.00 & 0.05 & 0.0826 & 0.02 & 0.06\\\hline
				0.25 & 200 & 0.1939 & 0.04 & 0.06 & 0.2029 & 0.03 & 0.03 & 0.1947 & 0.03 & 0.03\\
				& 500 & 0.1203 & 0.00 & 0.02 & 0.1305 & 0.00 & 0.02 & 0.1278 & 0.00 & 0.02\\
				& 1000 & 0.0941 & 0.01 & 0.03 & 0.0890 & 0.00 & 0.03 & 0.0932 & 0.00 & 0.05\\\hline
				0.5 & 200 & 0.2125 & 0.01 & 0.05 & 0.2272 & 0.00 & 0.04 & 0.2003 & 0.01 & 0.04\\
				& 500 & 0.1395 & 0.00 & 0.04 & 0.1352 & 0.00 & 0.04 & 0.1338 & 0.00 & 0.02\\
				& 1000 & 0.1004 & 0.01 & 0.06 & 0.1096 & 0.00 & 0.06 & 0.1104 & 0.01 & 0.05\\
				\hline\hline
				& & \multicolumn{3}{c|}{$\theta_{10}(t)=0$}\\
				\hline
				& & & \multicolumn{1}{c}{${\rm{H}}_{0{\rm{C}}}$} & \multicolumn{1}{c|}{${\rm{H}}_{0{\rm{Z}}}$}\\
				\multicolumn{1}{c}{$\delta$} & \multicolumn{1}{c|}{$n$} & \multicolumn{1}{c}{MSE} & \multicolumn{1}{c}{E1} & \multicolumn{1}{c|}{E1}\\
				\hline\hline
				0.1 & 200 & 0.1625 & 0.02 & 0.04\\
				& 500 & 0.1094 & 0.00 & 0.02\\
				& 1000 & 0.6732 & 0.00 & 0.01\\\hline
				0.25 & 200 & 0.1924 & 0.03 & 0.04\\
				& 500 & 0.1150 & 0.00 & 0.01\\
				& 1000 & 0.0890 & 0.00 & 0.01\\\hline
				0.5 & 200 & 0.1848 & 0.01 & 0.02\\
				& 500 & 0.1162 & 0.00 & 0.01\\
				& 1000 & 0.1002 & 0.00 & 0.05\\
				\hline
			\end{tabular}
		}
	\end{center}
\end{table}

The second model setting focuses on the case where $p$ is larger than in the first model setting. We set $p=10$ and $q=1$ and defined the coefficient functions $\theta_j(\cdot),\ j=1, 2, \ldots, 10$ as
\begin{equation*}
	\theta_j(t)=\begin{cases}
		1, & j=1,\\
		\cos(2t), & j=2,\\
		0, & j=3, 4, \ldots, 10,
	\end{cases}
\end{equation*}
where the intercept term $\theta_0(t)$ was not considered. The kernel function was the Epanechnikov kernel, and the tuning parameters were selected in the same manner as in the first model setting. Table \ref{tab2} shows the calculated MSEs and empirical probabilities of error for each coefficient function $\theta_j(\cdot)$ when $\delta=0.1, 0.25, 0.5$ and $n=200, 500, 1000$. The accuracy of the estimator and test statistic improved as $n$ increased, with no significant deterioration compared to the first model setting with $p=3$, indicating that the proposed model can avoid the curse of dimensionality even when the dimension $p$ is large. Figure \ref{fig1} shows the results of the estimation. The two dotted lines are plots of the 5th and 95th largest estimates of the $M=100$ estimates at each point $t\in[0, 1]$. The average estimates (dashed line) resembled the true value (solid line). 

\begin{table}[t]
	\begin{center}
		\caption{\footnotesize{Results of estimation and hypothesis testing in the third model setting. E1 and E2 represent the empirical probabilities of Type I and Type II errors, respectively. Since $\theta_0({\bf{t}})=2$ is a nonzero constant, the null hypothesis ${\rm{H}}_{0{\rm{C}}}$ is true, whereas the null hypothesis ${\rm{H}}_{0{\rm{Z}}}$ is false. Accordingly, we provide E1 for ${\rm{H}}_{0{\rm{C}}}$ and E2 for ${\rm{H}}_{0{\rm{Z}}}$. Meanwhile, since $\theta_1({\bf{t}})=-\exp(-10\|{\bf{t}}-(0.5, 0.5)^\top\|^2)$ is not constant, the null hypotheses ${\rm{H}}_{0{\rm{C}}}$ and ${\rm{H}}_{0{\rm{Z}}}$ are false and thus E2 is given for both tests. Similarly, $\theta_2({\bf{t}})=0$ indicates that the null hypotheses ${\rm{H}}_{0{\rm{C}}}$ and ${\rm{H}}_{0{\rm{Z}}}$ are true. Thus, for $\theta_2({\bf{t}})=0$, E1 is provided for both tests.\\}}\label{tab3}
		\begin{tabular}{@{}rr|rrr|rrr|rrr@{}}
			\hline
			& & \multicolumn{3}{c|}{$\theta_0({\bf{t}})=2$} & \multicolumn{3}{l|}{$\theta_1({\bf{t}})=-\exp(-10$} & \multicolumn{3}{c}{$\theta_2({\bf{t}})=0$}\\
			& & & & & \multicolumn{3}{r|}{$\times\|{\bf{t}}-(0.5, 0.5)^\top\|^2)$} & \\
			\hline
			& & & \multicolumn{1}{c}{${\rm{H}}_{0{\rm{C}}}$} & \multicolumn{1}{c|}{${\rm{H}}_{0{\rm{Z}}}$} & & \multicolumn{1}{c}{${\rm{H}}_{0{\rm{C}}}$} & \multicolumn{1}{c|}{${\rm{H}}_{0{\rm{Z}}}$} & & \multicolumn{1}{c}{${\rm{H}}_{0{\rm{C}}}$} & \multicolumn{1}{c}{${\rm{H}}_{0{\rm{Z}}}$}\\
			\multicolumn{1}{c}{$\delta$} & \multicolumn{1}{c|}{$n$} & \multicolumn{1}{c}{MSE} & \multicolumn{1}{c}{E1} & \multicolumn{1}{c|}{E2} & \multicolumn{1}{c}{MSE} & \multicolumn{1}{c}{E2} & \multicolumn{1}{c|}{E2} & \multicolumn{1}{c}{MSE} & \multicolumn{1}{c}{E1} & \multicolumn{1}{c}{E1}\\
			\hline\hline
			0.1 & 3000 & 0.9334 & 0.00 & 0.00 & 0.1475 & 0.16 & 0.00 & 0.0865 & 0.00 & 0.01\\
			& 5000 & 0.9378 & 0.00 & 0.00 & 0.1222 & 0.03 & 0.00 & 0.0695 & 0.00 & 0.00\\\hline
			0.25 & 3000 & 0.9033 & 0.00 & 0.00 & 0.1439 & 0.19 & 0.00 & 0.0881 & 0.01 & 0.00\\
			& 5000 & 0.9191 & 0.01 & 0.00 & 0.1163 & 0.02 & 0.00 & 0.0719 & 0.01 & 0.01\\\hline
			0.5 & 3000 & 0.8810 & 0.01 & 0.00 & 0.1408 & 0.15 & 0.00 & 0.0907 & 0.00 & 0.00\\
			& 5000 & 0.8827 & 0.00 & 0.00 & 0.1171 & 0.02 & 0.00 & 0.0704 & 0.00 & 0.00 \\
			\hline
		\end{tabular}
	\end{center}
\end{table}

In the third model setting, we set $p=2$ and $q=2$ and defined the coefficient functions $\theta_j(\cdot),\ j=0, 1, 2$ as
\begin{equation*}
	\theta_j({\bf{t}})=\begin{cases}
		2, & j=0,\\
		-\exp(-10\|{\bf{t}}-(0.5, 0.5)^\top\|^2), & j=1,\\
		0, & j=2.
	\end{cases}
\end{equation*}
We employed the kernel function of the Epanechnikov type as follows:
\begin{equation*}
	K({\bf{u}})=\frac{2}{\pi}(1-\|{\bf{u}}\|^2)I(\|{\bf{u}}\|\leq1).
\end{equation*}
The tuning parameters were selected in the same manner as in the first model setting. Table \ref{tab3} shows the calculated MSEs and empirical probabilities of error for each coefficient function $\theta_j(\cdot)$ when $\delta=0.1, 0.25, 0.5$ and $n=3000, 5000$. As with the first and second settings, the accuracy of the estimator and test statistic improved as $n$ increased.

We note that Tables \ref{tab1}-\ref{tab3} show the results of the hypothesis tests when the tuning parameters are automatically selected based on each dataset. As a result, the empirical probability of the Type I error tended to be smaller than the given significance level $\alpha=0.05$ in many settings, and thus the results seem to be conservative. However, although ad hoc selection of the tuning parameters will yield more reasonable results about the Type I error, such tuning parameters cannot be determined in real data analysis.

In this study, we used the plug-in test, but the bootstrap method may also be useful for testing. However, as described in de Haan and Zhou (2022), it is known that bootstrap methods do not always work efficiently in the context of extreme value theory. In addition, the effectiveness of bootstrapping in extreme value theory has only been partially revealed. Thus, at the present stage, the bootstrap test in our model is posited as future research.

\section{Application}\label{sec5}

\begin{figure}[t!]
	\centering
	\includegraphics[keepaspectratio, width=9.5cm,angle=270]{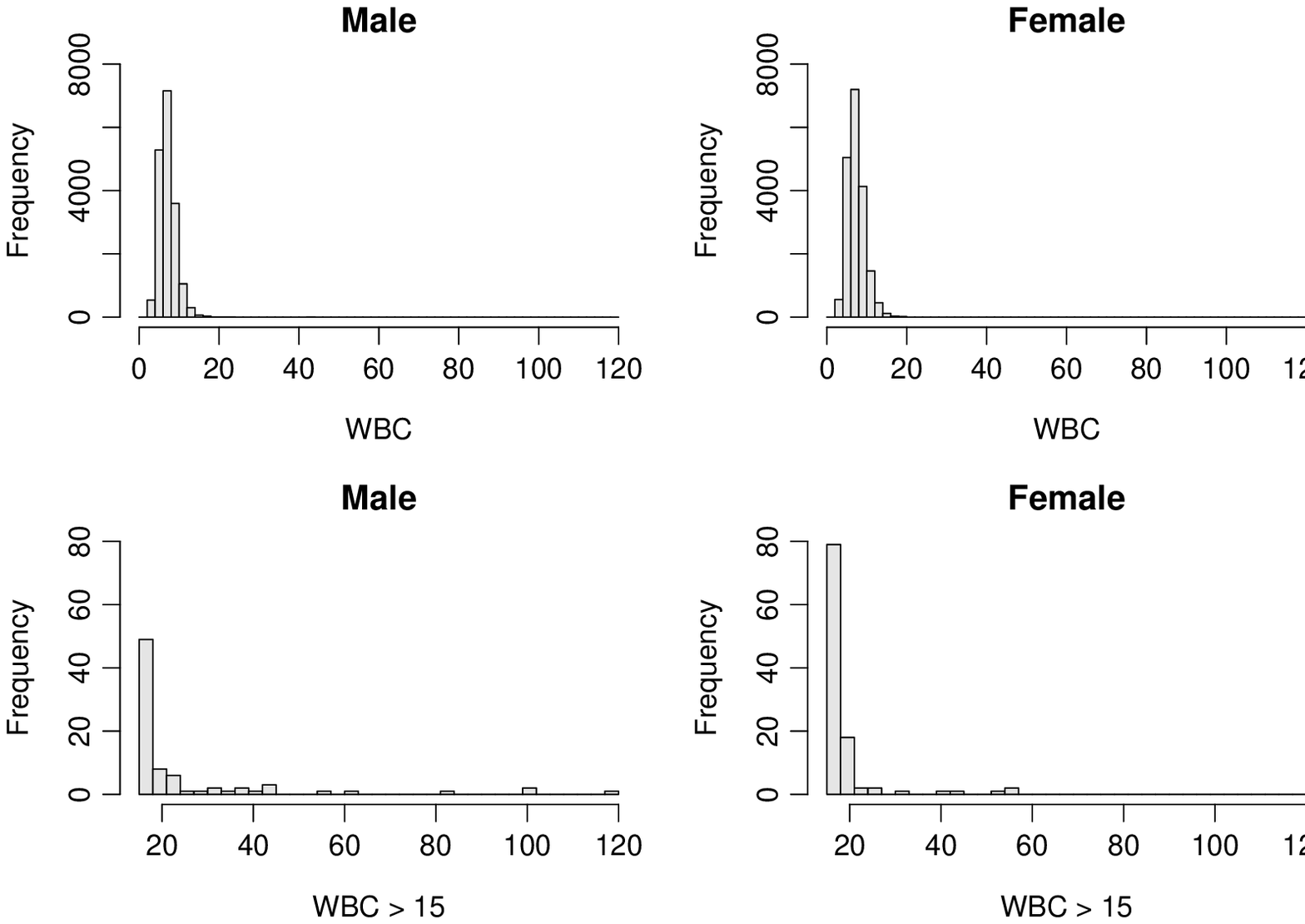}
	\caption{The histograms of the response $Y$ for male (top left panel) and female (top right panel): The bottom two panels show the histograms of the response $Y$ greater than 15 for male (bottom left panel) and female (bottom right panel).}
	\label{fig2}
\end{figure}

\begin{figure}[t!]
	\centering
	\includegraphics[keepaspectratio, width=9cm,angle=270]{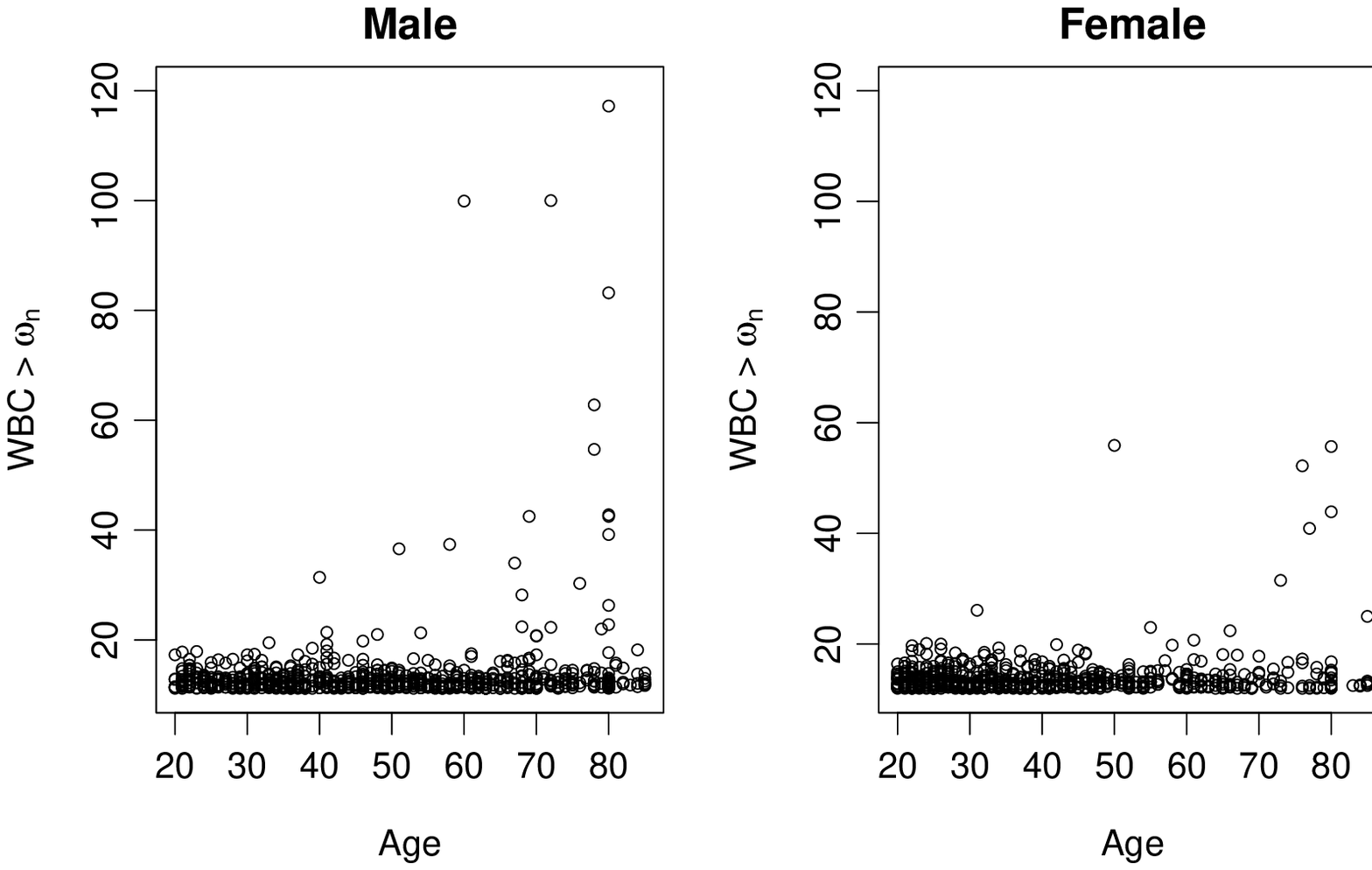}
	\caption{The time series plots of $Y$ for male (left panel) and female (right panel), where $Y$ exceeds the threshold $\omega_n$.}
	\label{fig3}
\end{figure}

In this section, we apply the proposed method to a real dataset on white blood cells. The dataset is available in Kaggle (\url{https://www.kaggle.com/amithasanshuvo/cardiac-data-nhanes}). White blood cells play a role in processing foreign substances such as bacteria and viruses that have invaded the body, and are a type of blood cell that is indispensable for maintaining the normal immune function of the human body. Therefore, if the white blood cell count is abnormal, diseases may be suspected. The top left and right panels of Figure \ref{fig2} show histograms of the white blood cell counts for $n=18047$ males and $n=19032$ females aged 20 to 85, respectively, and the bottom two panels show histograms for those over 15 ($\times 10^3/\mu$L). We can judge whether the tails of these distributions have a positive extreme value index by comparing them to the normal distribution with a zero extreme value index. In many extreme value studies, kurtosis is often used. The sample kurtosis was about $403.8$ for males and about $38.3$ for females, indicating that the right tails of these distributions are heavy. In addition, Figure \ref{fig3} shows plots of the subject's age and white blood cell count, suggesting that the number of abnormal white blood cell counts tends to increase with age.

The dataset also contains percentages by type: neutrophils, eosinophils, basophils, monocytes, and lymphocytes. White blood cell differentiation is a clinical test that identifies the types of white blood cells that cause an abnormal white blood cell count. These five types have different immune functions and can help detect certain diseases. The sample averages were about 58.02, 3.10, 0.69, 8.39, and 29.84\% for males and about 58.70, 2.57, 0.71, 7.47, and 30.59\% for females, respectively. Neutrophils and lymphocytes comprised the majority of white blood cells, and the correlation coefficient calculated from the transformed observations, as described below, was approximately $-0.93$ for males and $-0.95$ for females. In other words, there was a strong negative correlation between the percentages of neutrophils and lymphocytes. In this analysis, we define the response $Y$ as the white blood cell count; the predictors $X_1$, $X_2$, $X_3$ and $X_4$ as the percentages of eosinophils, basophils, monocytes and lymphocytes in the white blood cells; and the predictor $T$ as age. We denote ${\bf{X}}=(X_1, X_2, X_3, X_4)^\top$.

\begin{figure}[t!]
	\centering
	\includegraphics[keepaspectratio, width=9cm,angle=270]{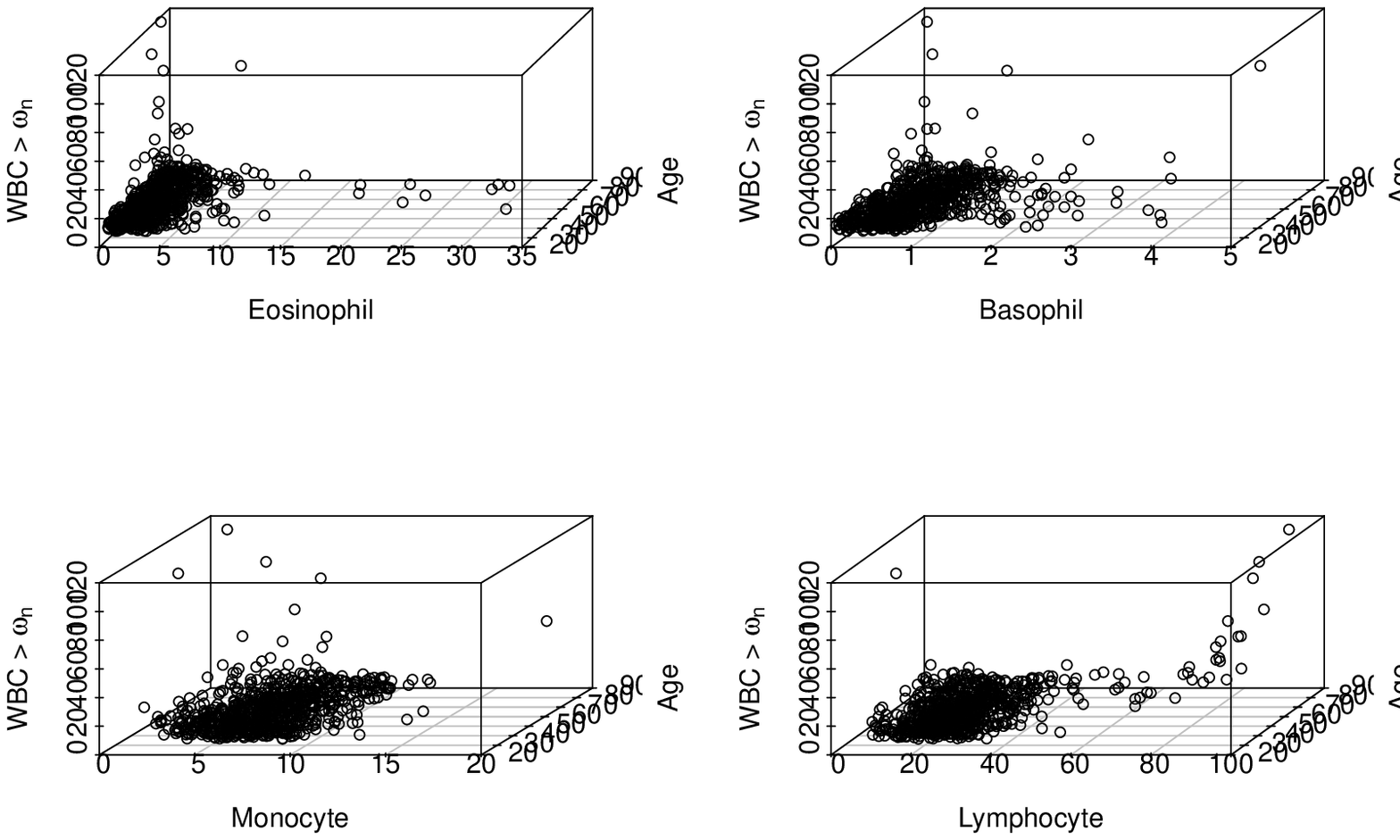}
	\caption{The three dimensional scatter plots of $(X_j, T, Y),\ j=1, 2, 3, 4$ with $Y>\omega_n$ for male. For the top left, top right, bottom left and bottom right panels, $X_j$ is the percentage of eosinophils, basophils, monocytes and lymphocytes in the white blood cells, respectively.}
	\label{fig4}
\end{figure}

\begin{figure}[t!]
	\centering
	\includegraphics[keepaspectratio, width=9cm,angle=270]{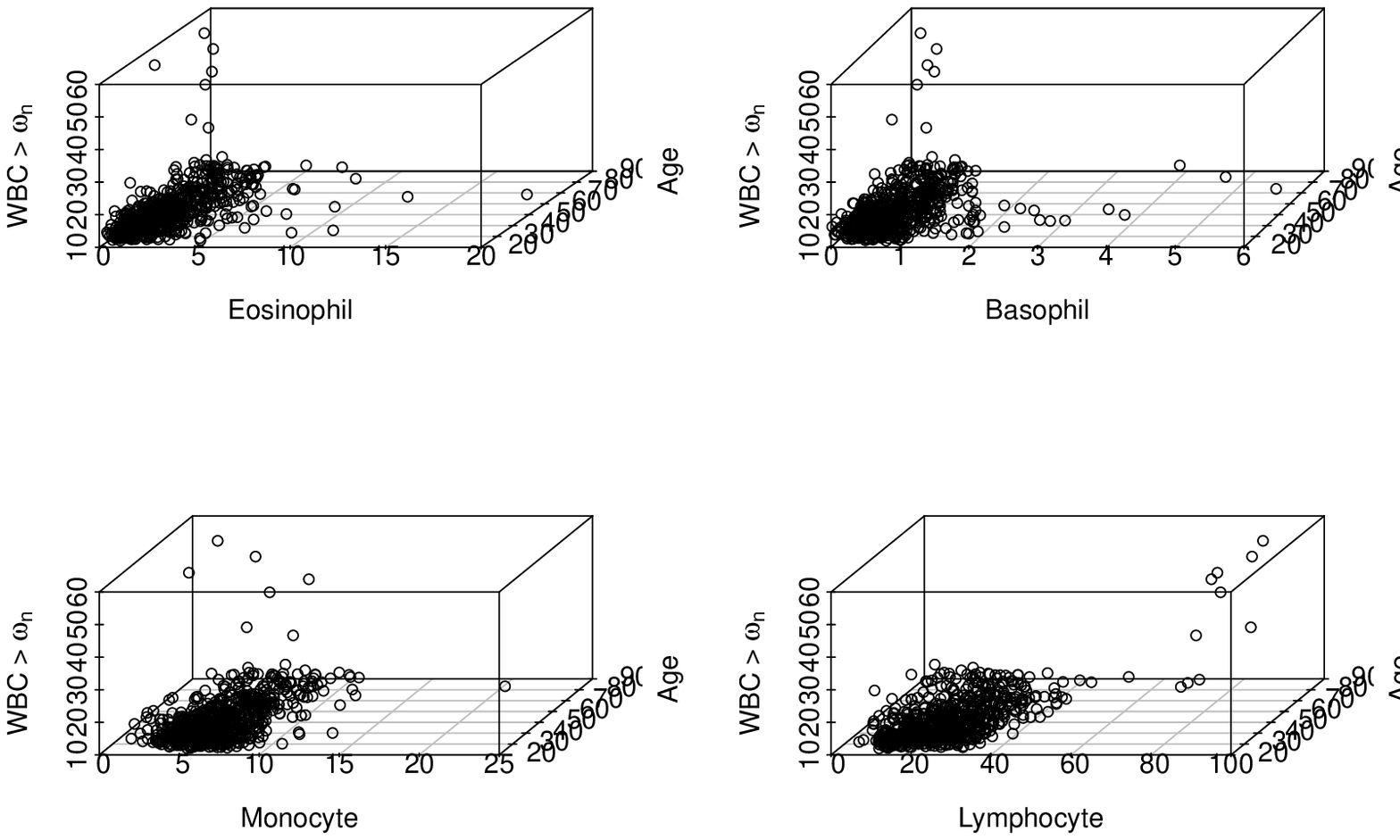}
	\caption{The three dimensional scatter plots of $(X_j, T, Y),\ j=1, 2, 3, 4$ with $Y>\omega_n$ for female. For the top left, top right, bottom left and bottom right panels, $X_j$ is the percentage of eosinophils, basophils, monocytes and lymphocytes in the white blood cells, respectively.}
	\label{fig5}
\end{figure}

Figures \ref{fig4} and \ref{fig5} show the three-dimensional scatter plots of each $(X_j, T, Y)$ for male and female, respectively. As shown in these figures, the predictors $X_1$, $X_2$, $X_3$ and $X_4$ had many outliers. However, excluding these outliers also excludes the extreme values of the response $Y$. Therefore, we apply the normal score transformation to $X_j$. That is, if $X_{ij}$ is the $R_{ij}$-th largest in the $j$th predictor sample $\{X_{ij}\}_{i=1}^n$, $X_{ij}$ is redefined as
\begin{equation*}
	X_{ij}=\Phi^{-1}((R_{ij}-3/8)/(n+1/4)),
\end{equation*}
where all observations are jittered by uniform noise before applying the normal score transformation. Consequently, the redefined predictors $X_1$, $X_2$, $X_3$, and $X_4$ are normally distributed. Wang and Tsai (2009) applied a similar transformation in their analysis of real data.

\begin{figure}[t!]
	\centering
	\includegraphics[keepaspectratio, width=19cm,angle=270]{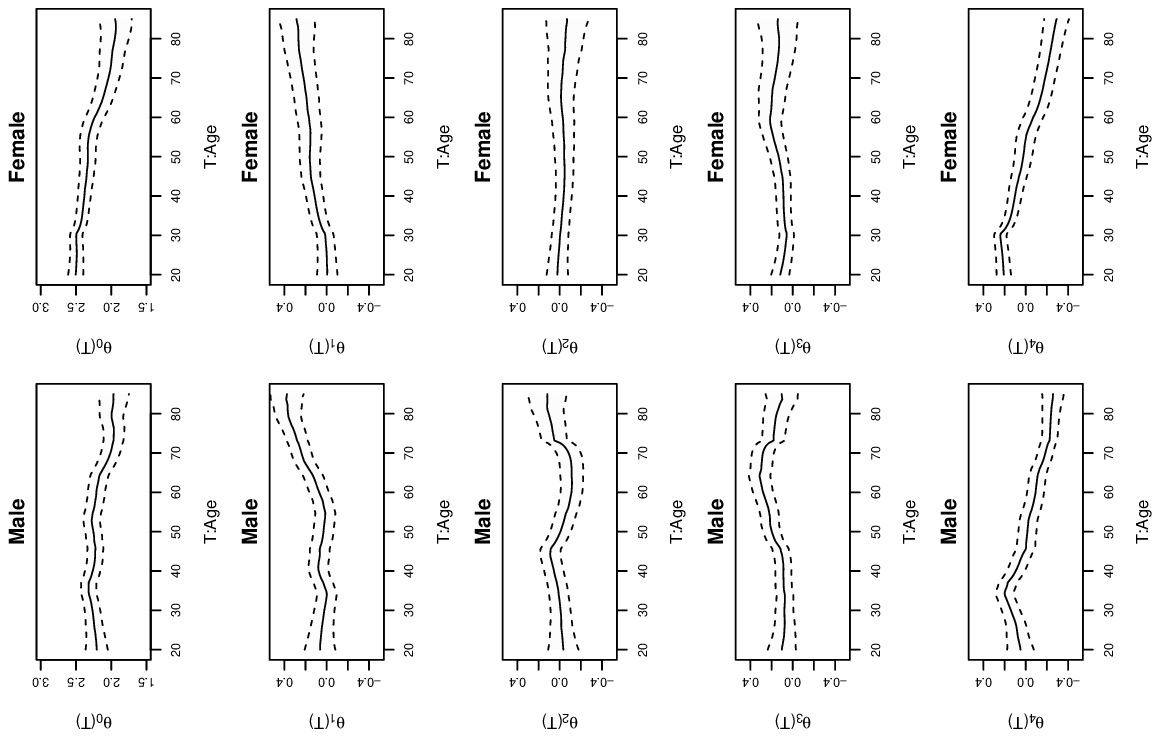}
	\caption{The estimated coefficient functions (solid line) and its 95\% confidence intervals (dashed lines) with bias ignored for male (first column) and female (second column).}
	\label{fig6}
\end{figure}

We assume that the conditional distribution function of $Y$ given $({\bf{X}}, T)=({\bf{x}}, t)$ satisfies
\begin{equation*}
	1-F(y\mid {\bf{x}}, t)=y^{-1/\gamma({\bf{x}}, t)}\mathcal{L}(y; {\bf{x}}, t),
\end{equation*}
where $\mathcal{L}(\cdot; {\bf{x}}, t)$ is a slowly varying function satisfying (\ref{2.1.2}), and
\begin{equation}
	\log\gamma({\bf{x}}, t)^{-1}=\theta_0(t)+\theta_1(t)x_1+\theta_2(t)x_2+\theta_3(t)x_3+\theta_4(t)x_4,\label{5.1}\tag{5.1}
\end{equation}
where ${\bf{x}}=(x_1, x_2, x_3, x_4)^\top\in\mathbb{R}^4$, and $\theta_j(t),\ j=0, 1, 2, 3, 4$ are unknown smooth functions of $t$. The aim of the analysis is to investigate the effect of $X_j$ on the extreme values of $Y$, where the effect of $X_j$ varies with $T$. To do this, we first estimate the unknown coefficient functions $\theta_j(\cdot),\ j=0, 1, 2, 3, 4$. Then, we select the threshold $\omega_n$ and bandwidth $h_n$ using the procedure described in Section \ref{subsec2.3}. We employ the Epanechikov kernel in the proposed estimator and set the pre-determined sample fraction to $n_0/n=0.030$ in $D=20$-fold cross-validation, where $n_0=\sum_{i=1}^nI(Y_i>\omega_n)$. We obtained the optimal tuning parameters as $(h_n, n_0/n)=(0.21, 0.042)$ for male and $(h_n, n_0/n)=(0.30, 0.036)$ for female. Figure \ref{fig6} shows the estimated coefficient functions $\widehat{\theta}_j(\cdot),\ j=0, 1, 2, 3, 4$ by the solid line and the following pointwise 95\% confidence intervals computed from the asymptotic normality of the proposed estimator by the dashed lines:
\begin{align*}
	&\left(\widehat{\theta}_j(t)-\Phi^{-1}(0.975)\nu^{1/2}[\widehat{n(t)\sigma_{nj}(t)}]^{-1/2},\right.\\
	&\hspace*{3cm}\left.\widehat{\theta}_j(t)+\Phi^{-1}(0.975)\nu^{1/2}[\widehat{n(t)\sigma_{nj}(t)}]^{-1/2}\right),
\end{align*}
where the bias is ignored, $n(t)\sigma_{nj}(t)$ defined in Section \ref{subsec4.1} is estimated based on (C.3), and $\nu=\int K(u)^2{\rm{d}}u$. For all coefficient functions, the trends were similar for male and female. The decreasing trend in the estimated intercept term $\widehat{\theta}_0(\cdot)$ indicates that the number of abnormal white blood cell counts tends to increase with age. Some of the estimated coefficient functions deviated from zero and vary with age.

\begin{table}[t!]
	\begin{center}
		\caption{\footnotesize{The results of the hypothesis testing. ${\rm{H}}_{0\rm{C}}$ is the null hypothesis that $\theta_j(t)$ is constant, and ${\rm{H}}_{0\rm{Z}}$ is the null hypothesis that $\theta_j(t)$ is zero. For the significance level $\alpha=0.05$, we reject the null hypothesis if $\widetilde{T}<-0.61$ or $\widetilde{T}>4.37$.\\}}\label{tab4}
		\begin{tabular}{@{}cccrcc@{}}
			\hline
			Coefficient & Null&\multirow{2}{*}{Gender} & \multicolumn{1}{c}{Test} & \multirow{2}{*}{p-value} & \multirow{2}{*}{Result}\\
			function & hypothesis & & \multicolumn{1}{c}{statistic $\widetilde{T}$} & &\\
			\hline\hline
			$\theta_0$ & ${\rm{H}}_{0\rm{C}}$ & Male & 3.25 & 0.07477 & \\
			& & Female & 7.31 & 0.00134 & Rejected\\
			& ${\rm{H}}_{0\rm{Z}}$ & Male & 72.63 & 0.00000 & Rejected\\
			& & Female & 85.45 & 0.00000 & Rejected\\\hline
			$\theta_1$ & ${\rm{H}}_{0\rm{C}}$ & Male & 5.15 & 0.01150 & Rejected\\
			& & Female & 4.25 & 0.02816 &\\
			& ${\rm{H}}_{0\rm{Z}}$ & Male & 8.79 & 0.00031 & Rejected\\
			& & Female & 4.75 & 0.01712 & Rejected\\\hline
			$\theta_2$ & ${\rm{H}}_{0\rm{C}}$ & Male & 2.00 & 0.23653 &\\
			& & Female & 0.48 & 0.28831 &\\
			& ${\rm{H}}_{0\rm{Z}}$ & Male & 2.09 & 0.21855 &\\
			& & Female & 0.71 & 0.62501 &\\\hline
			$\theta_3$ & ${\rm{H}}_{0\rm{C}}$ & Male & 3.33 & 0.06894 &\\
			& & Female & 2.26 & 0.18850 &\\
			& ${\rm{H}}_{0\rm{Z}}$ & Male & 8.97 & 0.00025 & Rejected\\
			& & Female & 5.16 & 0.01143 & Rejected\\\hline
			$\theta_4$ & ${\rm{H}}_{0\rm{C}}$ & Male & 7.94 & 0.00071 & Rejected\\
			& & Female & 11.97 & 0.00001 & Rejected\\
			& ${\rm{H}}_{0\rm{Z}}$ & Male & 9.00 & 0.00025 & Rejected\\
			& & Female & 11.42 & 0.00002 & Rejected\\
			\hline
		\end{tabular}
	\end{center}
\end{table}

Table \ref{tab4} presents the results of the statistical hypothesis tests for sparsity and constancy, as defined in Section \ref{subsec4.1}. For the significance level $\alpha=0.05$, we reject the null hypothesis if $\widetilde{T}<-0.61$ or $\widetilde{T}>4.37$. The null hypothesis ${\rm{H}}_{0\rm{Z}}$ for sparsity was rejected for all coefficient functions, except $\theta_2(\cdot)$ for both male and female. In addition, the null hypothesis ${\rm{H}}_{0\rm{C}}$ for constancy was rejected for $\theta_1(\cdot)$ and $\theta_4(\cdot)$ for male and $\theta_0(\cdot)$ and $\theta_4(\cdot)$ for female. Remarkably, eosinophils and monocytes, which represented a small percentage of white blood cells, were associated with abnormal white blood cell counts. 

\begin{figure}[t!]
	\centering
	\includegraphics[keepaspectratio, width=7cm,angle=270]{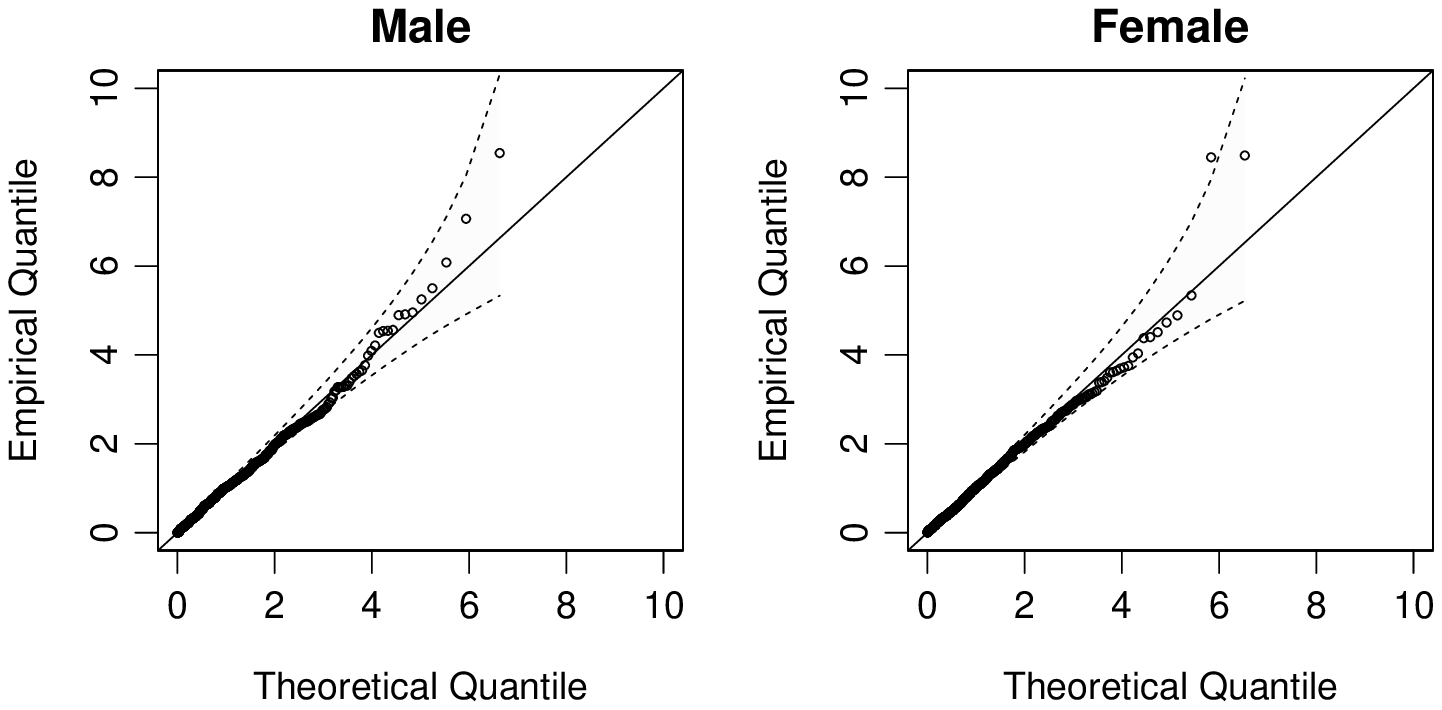}
	\caption{The Q-Q plots for the proposed model for male (left panel) and female (right panel).}
	\label{fig7}
\end{figure}

\begin{figure}[t!]
	\centering
	\includegraphics[keepaspectratio, width=7cm,angle=270]{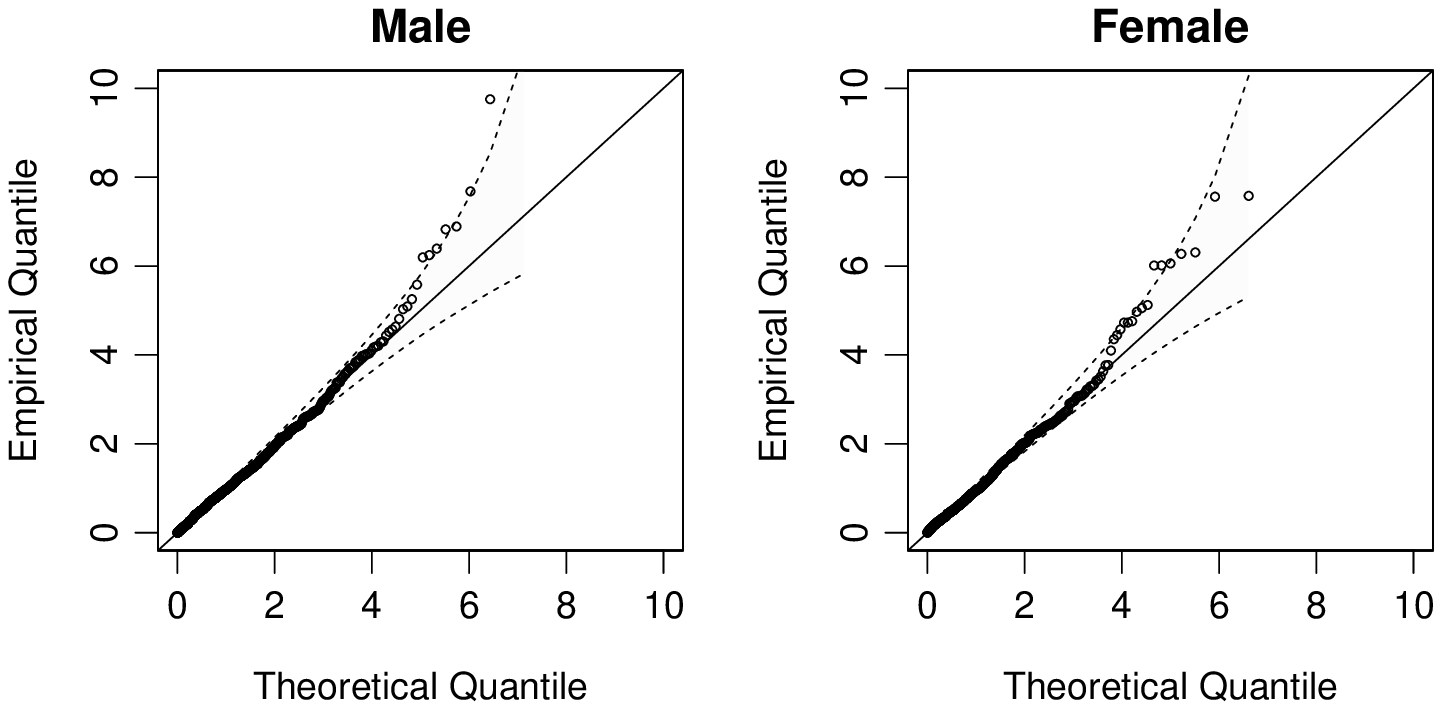}
	\caption{The Q-Q plots for the linear model proposed by Wang and Tsai (2009) for male (left panel) and female (right panel).}
	\label{fig8}
\end{figure}

We evaluate the goodness of fit of the model using the Q-Q plot (quantile-quantile plot). We regard $\{\exp((1, {\bf{X}}_i^\top)\widehat{\bm{\theta}}(T_i))\log(Y_i/\omega_n): Y_i>\omega_n,\ i=1, \ldots, n\}$ as a random sample from the standard exponential distribution. Figure \ref{fig7} shows plots of these empirical and theoretical quantiles. The two dashed lines show the pointwise 95\% confidence intervals computed in the simulations. We can infer that the better the plots are aligned on a straight line, the better the model fits the data. Most of the plots were within the 95\% confidence interval and the goodness of fit of the model did not seem to be bad. In contrast, Figure \ref{fig8} shows the plots for the linear model proposed by Wang and Tsai (2009), where the predictors are defined as $T$ scaled on $[0,1]$, $X_1$, $X_2$, $X_3$ and $X_4$. In this case, many plots were outside the 95\% confidence interval and deviated significantly from the straight line, indicating that our model fits the data better.

Finally, because the null hypotheses ${\rm{H}}_{0\rm{Z}}$ and ${\rm{H}}_{0\rm{C}}$ were not rejected for $\theta_2(\cdot)$ and $\theta_3(\cdot)$, we adopt a simpler model. We consider the model
\begin{equation}
	\log\gamma({\bf{x}}, t)^{-1}=\theta_0(t)+\theta_1(t)x_1+\theta_3(t)x_3+\theta_4(t)x_4,\label{5.2}\tag{5.2}
\end{equation}
which assumes the sparsity of $X_2$. For the model (\ref{5.1}), the discrepancy measure value described in Section \ref{subsec2.3} was approximately $4.322\times10^{-4}$ for males and $3.015\times10^{-4}$ for females. Meanwhile, for the model (\ref{5.2}), the discrepancy measure value was approximately $3.730\times10^{-4}$ for males and $3.017\times10^{-4}$ for females, where $(h_n, n_0/n)=(0.22, 0.042)$ for males and $(h_n, n_0/n)=(0.30, 0.036)$ for females. The discrepancy measure values for females were not very different between the two models, but the discrepancy measure value for males was smaller in the model (\ref{5.2}) than in the model (\ref{5.1}). Moreover, we consider the model
\begin{equation}
	\log\gamma({\bf{x}}, t)^{-1}=\theta_0(t)+\theta_1(t)x_1+\widehat{\theta}_3x_3+\theta_4(t)x_4,\label{5.3}\tag{5.3}
\end{equation}
where $\widehat{\theta}_3$ is the average of the estimates $\{\widehat{\theta}_3(t_l)\}_{l=1}^L$ obtained in model (\ref{5.1}), which is a known constant. For the model (\ref{5.3}), the discrepancy measure value was approximately $3.628\times10^{-4}$ for males and $3.104\times10^{-4}$ for females, where $(h_n, n_0/n)=(0.19, 0.042)$ for males and $(h_n, n_0/n)=(0.30, 0.036)$ for females. The discrepancy measure value for males was smaller in the model (\ref{5.3}) than in the model (\ref{5.1}). Therefore, from the point of view of the discrepancy measure, the data structure may be explained by a simpler model.

\section*{Appendix}

In this appendix, we prove Theorems \ref{thm1}-\ref{thm3} for ${\bf{t}}={\bf{t}}_0\in\mathbb{R}^q$. For convenience, the intercept term $\theta_0(\cdot)$ is not considered.

\bigskip

\begin{proof}[Proof of Theorem~{\upshape\ref{thm1}}]
	$\dot{\bm{L}}_n({\bm{\theta}}({\bf{t}}_0))$ can be regarded as the sum of independent and identically distributed random variables. To apply the Central Limit Theorem, we show $E[[n_0({\bf{t}}_0){\bm{\Sigma}}_n({\bf{t}}_0)]^{-1/2}(\dot{\bm{L}}_n({\bm{\theta}}({\bf{t}}_0))+n_0({\bf{t}}_0){\bm{\Sigma}}_n({\bf{t}}_0)\sum_{l=1}^2{\bf{\Lambda}}_n^{(l)}({\bf{t}}_0))]\to-{\bf{b}}({\bf{t}}_0)$ as $n\to\infty$ in the first step and ${\rm{var}}[[n_0({\bf{t}}_0){\bm{\Sigma}}_n({\bf{t}}_0)]^{-1/2}(\dot{\bm{L}}_n({\bm{\theta}}({\bf{t}}_0))+n_0({\bf{t}}_0){\bm{\Sigma}}_n({\bf{t}}_0)\sum_{l=1}^2{\bf{\Lambda}}_n^{(l)}({\bf{t}}_0))]={\rm{var}}[[n_0({\bf{t}}_0){\bm{\Sigma}}_n({\bf{t}}_0)]^{-1/2}\dot{\bm{L}}_n({\bm{\theta}}({\bf{t}}_0))]\to\nu{\bf{I}}_p$ as $n\to\infty$ in the second step, where ``var'' denotes the variance-covariance matrix.
	
	\begin{step}\label{step1}
		We can write $\dot{\bm{L}}_n({\bm{\theta}}({\bf{t}}_0))$ as
		\begin{align*}
			\dot{\bm{L}}_n({\bm{\theta}}({\bf{t}}_0))&=\sum_{i=1}^n{\bf{X}}_i\left\{\exp({\bf{X}}_i^\top{\bm{\theta}}({\bf{t}}_0))\log\frac{Y_i}{\omega_n}-1\right\}I(Y_i>\omega_n)K({\bf{H}}_n^{-1}({\bf{t}}_0-{\bf{T}}_i))\\
			&=\sum_{i=1}^n{\bf{X}}_i\left\{\exp({\bf{X}}_i^\top{\bm{\theta}}({\bf{T}}_i))\log\frac{Y_i}{\omega_n}-1\right\}I(Y_i>\omega_n)K({\bf{H}}_n^{-1}({\bf{t}}_0-{\bf{T}}_i))\\
			&\quad-\sum_{i=1}^n{\bf{X}}_i\left\{\exp({\bf{X}}_i^\top{\bm{\theta}}({\bf{T}}_i))-\exp({\bf{X}}_i^\top{\bm{\theta}}({\bf{t}}_0))\right\}\\
			&\quad\quad\times\log\frac{Y_i}{\omega_n}I(Y_i>\omega_n)K({\bf{H}}_n^{-1}({\bf{t}}_0-{\bf{T}}_i))\\
			&=:{\bf{A}}_n^{(1)}-{\bf{A}}_n^{(2)}.
		\end{align*}
		By the Taylor expansion and the condition (C.1), we have
		\begin{align*}
			&E[[n_0({\bf{t}}_0){\bf{\Sigma}}_n({\bf{t}}_0)]^{-1/2}{\bf{A}}_n^{(1)}]\\
			&=n[n_0({\bf{t}}_0){\bf{\Sigma}}_n({\bf{t}}_0)]^{-1/2}\\
			&\quad\times\int_{\omega_n}^\infty\int_{\mathbb{R}^p}\int_{\mathbb{R}^q}{\bf{x}}\left(\gamma({\bf{x}}, {\bf{t}})^{-1}\log\frac{y}{\omega_n}-1\right)K({\bf{H}}_n^{-1}({\bf{t}}_0-{\bf{t}}))f(y, {\bf{x}}, {\bf{t}}){\rm{d}}y{\rm{d}}{\bf{x}}{\rm{d}}{\bf{t}}
		\end{align*}
		\begin{align*}
			&=n\det({\bf{H}}_n)[n_0({\bf{t}}_0){\bf{\Sigma}}_n({\bf{t}}_0)]^{-1/2}\\
			&\quad\times\int_{\omega_n}^\infty\int_{\mathbb{R}^p}\int_{\mathbb{R}^q}{\bf{x}}\left(\gamma({\bf{x}}, {\bf{t}}_0+{\bf{H}}_n{\bf{s}})^{-1}\log\frac{y}{\omega_n}-1\right)K({\bf{s}})f(y, {\bf{x}}, {\bf{t}}_0+{\bf{H}}_n{\bf{s}}){\rm{d}}y{\rm{d}}{\bf{x}}{\rm{d}}{\bf{s}}\\
			&=n({\bf{t}}_0)[n_0({\bf{t}}_0){\bf{\Sigma}}_n({\bf{t}}_0)]^{-1/2}\int_{\omega_n}^\infty\int_{\mathbb{R}^p}{\bf{x}}\gamma({\bf{x}}, {\bf{t}}_0)^{-1}\log\frac{y}{\omega_n}\frac{f(y, {\bf{x}}, {\bf{t}}_0)}{f_{\bf{T}}({\bf{t}}_0)}{\rm{d}}y{\rm{d}}{\bf{x}}\\
			&\quad-n({\bf{t}}_0)[n_0({\bf{t}}_0){\bf{\Sigma}}_n({\bf{t}}_0)]^{-1/2}\int_{\omega_n}^\infty\int_{\mathbb{R}^p}{\bf{x}}\frac{f(y, {\bf{x}}, {\bf{t}}_0)}{f_{\bf{T}}({\bf{t}}_0)}{\rm{d}}y{\rm{d}}{\bf{x}}+o(1)\\
			&=:{\bf{B}}_n^{(1)}-{\bf{B}}_n^{(2)}+o(1).
		\end{align*}
		From the condition (C.4) and model assumptions (\ref{2.1.1}) and (\ref{2.1.2}), we have
		\begin{align*}
			{\bf{B}}_n^{(1)}&=n({\bf{t}}_0)[n_0({\bf{t}}_0){\bf{\Sigma}}_n({\bf{t}}_0)]^{-1/2}\int_{\omega_n}^\infty\int_{\mathbb{R}^p}{\bf{x}}\gamma({\bf{x}}, {\bf{t}}_0)^{-1}\log\frac{y}{\omega_n}\frac{f(y, {\bf{x}}, {\bf{t}}_0)}{f_{\bf{T}}({\bf{t}}_0)}{\rm{d}}y{\rm{d}}{\bf{x}}\\
			&=n({\bf{t}}_0)[n_0({\bf{t}}_0){\bf{\Sigma}}_n({\bf{t}}_0)]^{-1/2}\\
			&\quad\times E\left[{\bf{X}}\gamma({\bf{X}}, {\bf{t}}_0)^{-1}E\left[\log\frac{Y}{\omega_n}I(Y>\omega_n)\mid {\bf{X}}={\bf{X}}, {\bf{T}}={\bf{t}}_0\right]\mid {\bf{T}}={\bf{t}}_0\right]\\
			&=n({\bf{t}}_0)[n_0({\bf{t}}_0){\bf{\Sigma}}_n({\bf{t}}_0)]^{-1/2}\\
			&\quad\times E\left[{\bf{X}}\gamma({\bf{X}}, {\bf{t}}_0)^{-1}\int_0^\infty(1-F(\omega_ne^s\mid {\bf{X}}, {\bf{t}}_0)){\rm{d}}s\mid {\bf{T}}={\bf{t}}_0\right]\\
			&=n({\bf{t}}_0)[n_0({\bf{t}}_0){\bf{\Sigma}}_n({\bf{t}}_0)]^{-1/2}\\
			&\quad\times\left\{E\left[{\bf{X}}c_0({\bf{X}}, {\bf{t}}_0)\omega_n^{-1/\gamma({\bf{X}}, {\bf{t}}_0)}\int_0^\infty\gamma({\bf{X}}, {\bf{t}}_0)^{-1}e^{-s/\gamma({\bf{X}}, {\bf{t}}_0)}{\rm{d}}s\mid {\bf{T}}={\bf{t}}_0\right]\right.\\
			&\quad\quad+E\biggl[{\bf{X}}c_1({\bf{X}}, {\bf{t}}_0)\omega_n^{-1/\gamma({\bf{X}}, {\bf{t}}_0)-\beta({\bf{X}}, {\bf{t}}_0)}\biggr.\\
			&\quad\quad\quad\left.\left.\times\int_0^\infty\gamma({\bf{X}}, {\bf{t}}_0)^{-1}e^{-s(1/\gamma({\bf{X}}, {\bf{t}}_0)+\beta({\bf{X}}, {\bf{t}}_0))}{\rm{d}}s\mid {\bf{T}}={\bf{t}}_0\right]\right\}\{1+o(1)\}\\
			&=n({\bf{t}}_0)[n_0({\bf{t}}_0){\bf{\Sigma}}_n({\bf{t}}_0)]^{-1/2}\biggl\{E\left[{\bf{X}}c_0({\bf{X}},{\bf{t}}_0)\omega_n^{-1/\gamma({\bf{X}}, {\bf{t}}_0)}\mid {\bf{T}}={\bf{t}}_0\right]\biggr.\\
			&\quad+\biggl.E\left[{\bf{X}}\frac{c_1({\bf{X}}, {\bf{t}}_0)}{1+\gamma({\bf{X}}, {\bf{t}}_0)\beta({\bf{X}}, {\bf{t}}_0)}\omega_n^{-1/\gamma({\bf{X}}, {\bf{t}}_0)-\beta({\bf{X}}, {\bf{t}}_0)}\mid {\bf{T}}={\bf{t}}_0\right]\biggr\}\{1+o(1)\}
		\end{align*}
		Analogously, we have
		\begin{align*}
			{\bf{B}}_n^{(2)}&=n({\bf{t}}_0)[n_0({\bf{t}}_0){\bf{\Sigma}}_n({\bf{t}}_0)]^{-1/2}\int_{\omega_n}^\infty\int_{\mathbb{R}^p}{\bf{x}}\frac{f(y, {\bf{x}}, {\bf{t}}_0)}{f_{\bf{T}}({\bf{t}}_0)}{\rm{d}}y{\rm{d}}{\bf{x}}\\
			&=n({\bf{t}}_0)[n_0({\bf{t}}_0){\bf{\Sigma}}_n({\bf{t}}_0)]^{-1/2}\left\{E\left[{\bf{X}}c_0({\bf{X}}, {\bf{t}}_0)\omega_n^{-1/\gamma({\bf{X}}, {\bf{t}}_0)}\mid {\bf{T}}={\bf{t}}_0\right]\right.\\
			&\quad\left.+E\left[{\bf{X}}c_1({\bf{X}}, {\bf{t}}_0)\omega_n^{-1/\gamma({\bf{X}}, {\bf{t}}_0)-\beta({\bf{X}}, {\bf{t}}_0)}\mid {\bf{T}}={\bf{t}}_0\right]\right\}\{1+o(1)\}.
		\end{align*}
		Under the condition (C.5), we have
		\begin{align*}
			&{\bf{B}}_n^{(1)}-{\bf{B}}_n^{(2)}\\
			&=-n({\bf{t}}_0)[n_0({\bf{t}}_0){\bf{\Sigma}}_n({\bf{t}}_0)]^{-1/2}\\
			&\quad\times E\left[{\bf{X}}c_1({\bf{X}}, {\bf{t}}_0)\frac{\gamma({\bf{X}}, {\bf{t}}_0)\beta({\bf{X}}, {\bf{t}}_0)}{1+\gamma({\bf{X}}, {\bf{t}}_0)\beta({\bf{X}}, {\bf{t}}_0)}\omega_n^{-1/\gamma({\bf{X}}, {\bf{t}}_0)-\beta({\bf{X}}, {\bf{t}}_0)}\mid {\bf{T}}={\bf{t}}_0\right]+o(1)\\
			&\to-{\bf{b}}({\bf{t}}_0)
		\end{align*}
		as $n\to\infty$. Therefore, we have $E[[n_0({\bf{t}}_0){\bf{\Sigma}}_n({\bf{t}}_0)]^{-1/2}{\bf{A}}_n^{(1)}]\to-{\bf{b}}({\bf{t}}_0)$ as $n\to\infty$. Using the second-order Taylor expansion, we have
		\begin{align*}
			&\exp({\bf{X}}_i^\top{\bm{\theta}}({\bf{t}}_0+{\bf{H}}_n{\bf{s}}))-\exp({\bf{X}}_i^\top{\bm{\theta}}({\bf{t}}_0))\\
			&=\biggl\{{\bf{s}}^\top{\bf{H}}_n\exp({\bf{X}}_i^\top{\bm{\theta}}({\bf{t}}_0))({\bm{\Delta}}_1({\bf{t}}_0), \ldots, {\bm{\Delta}}_q({\bf{t}}_0))^\top{\bf{X}}_i\biggr.\\
			&\quad+\frac{1}{2}{\bf{s}}^\top{\bf{H}}_n\exp({\bf{X}}_i^\top{\bm{\theta}}({\bf{t}}_0))({\bm{\Delta}}_1({\bf{t}}_0), \ldots, {\bm{\Delta}}_q({\bf{t}}_0))^\top{\bf{X}}_i{\bf{X}}_i^\top({\bm{\Delta}}_1({\bf{t}}_0), \ldots, {\bm{\Delta}}_q({\bf{t}}_0)){\bf{H}}_n{\bf{s}}\\
			&\quad+\left.\frac{1}{2}{\bf{s}}^\top{\bf{H}}_n\exp({\bf{X}}_i^\top{\bm{\theta}}({\bf{t}}_0))[{\bf{X}}_i^\top{\bm{\Delta}}_{k_1k_2}({\bf{t}}_0)]_{q\times q}{\bf{H}}_n{\bf{s}}\right\}\{1+o(1)\}\\
			&=:\biggl\{\exp({\bf{X}}_i^\top{\bm{\theta}}({\bf{t}}_0)){\bf{s}}^\top {\bf{H}}_n{\bm{\lambda}}^{(0)}({\bf{X}}_i, {\bf{t}}_0)+\frac{1}{2}\exp({\bf{X}}_i^\top{\bm{\theta}}({\bf{t}}_0)){\bf{s}}^\top{\bf{H}}_n{\bm{\lambda}}^{(1)}({\bf{X}}_i, {\bf{t}}_0){\bf{H}}_n{\bf{s}}\biggr.\\
			&\hspace*{0.5cm}+\left.\frac{1}{2}\exp({\bf{X}}_i^\top{\bm{\theta}}({\bf{t}}_0)){\bf{s}}^\top{\bf{H}}_n{\bm{\lambda}}^{(2)}({\bf{X}}_i, {\bf{t}}_0){\bf{H}}_n{\bf{s}}\right\}\{1+o(1)\}.
		\end{align*}
		Therefore, by the Taylor expansion and condition (C.1), we have
		\begin{align}
			&E[[n_0({\bf{t}}_0){\bf{\Sigma}}_n({\bf{t}}_0)]^{-1/2}{\bf{A}}_n^{(2)}]\nonumber\\
			&=n[n_0({\bf{t}}_0){\bf{\Sigma}}_n({\bf{t}}_0)]^{-1/2}\int_{\omega_n}^\infty\int_{\mathbb{R}^p}\int_{\mathbb{R}^q}{\bf{x}}\left\{\exp({\bf{x}}^\top{\bm{\theta}}({\bf{t}}))-\exp({\bf{x}}^\top{\bm{\theta}}({\bf{t}}_0))\right\}\nonumber\\
			&\quad\times\log\frac{y}{\omega_n}K({\bf{H}}_n^{-1}({\bf{t}}_0-{\bf{t}}))f(y, {\bf{x}}, {\bf{t}}){\rm{d}}y{\rm{d}}{\bf{x}}{\rm{d}}{\bf{t}}\nonumber\\
			&=n\det({\bf{H}}_n)[n_0({\bf{t}}_0){\bf{\Sigma}}_n({\bf{t}}_0)]^{-1/2}\nonumber\\
			&\quad\times\int_{\omega_n}^\infty\int_{\mathbb{R}^p}\int_{\mathbb{R}^q}{\bf{x}}\left\{\exp({\bf{x}}^\top{\bm{\theta}}({\bf{t}}_0+{\bf{H}}_n{\bf{s}}))-\exp({\bf{x}}^\top{\bm{\theta}}({\bf{t}}_0))\right\}\nonumber\\
			&\quad\quad\times\log\frac{y}{\omega_n}K({\bf{s}})f(y, {\bf{x}}, {\bf{t}}_0+{\bf{H}}_n{\bf{s}}){\rm{d}}y{\rm{d}}{\bf{x}}{\rm{d}}{\bf{s}}\nonumber\\
			&=\frac{1}{2}n\det({\bf{H}}_n)[n_0({\bf{t}}_0){\bf{\Sigma}}_n({\bf{t}}_0)]^{-1/2}\nonumber\\
			&\quad\times\int_{\omega_n}^\infty\int_{\mathbb{R}^p}\int_{\mathbb{R}^q}{\bf{x}}\sum_{l=1}^2\left\{\gamma({\bf{x}}, {\bf{t}}_0)^{-1}{\bf{s}}^\top{\bf{H}}_n{\bm{\lambda}}^{(l)}({\bf{x}}, {\bf{t}}_0){\bf{H}}_n{\bf{s}}\right\}\nonumber\\
			&\quad\quad\times\log\frac{y}{\omega_n}K({\bf{s}})f(y, {\bf{x}}, {\bf{t}}_0){\rm{d}}y{\rm{d}}{\bf{x}}{\rm{d}}{\bf{s}}\{1+o(1)\}.\label{A.1}\tag{A.1}
		\end{align}
		Because the conditional distribution of $\gamma({\bf{X}}, {\bf{t}}_0)^{-1}\log(Y/\omega_n)$ given $({\bf{X}}, {\bf{T}})=({\bf{x}}, {\bf{t}}_0)$ and $Y>\omega_n$ is approximately a standard exponential, we have
		\begin{equation*}
			\int_{\omega_n}^\infty\gamma({\bf{x}},{\bf{t}}_0)^{-1}\log\frac{y}{\omega_n}\frac{f(y, {\bf{x}}, {\bf{t}}_0)/f_{({\bf{X}}, {\bf{T}})}({\bf{x}}, {\bf{t}}_0)}{1-F(\omega_n\mid {\bf{x}}, {\bf{t}}_0)}{\rm{d}}y\approx1,
		\end{equation*}
		where $f_{({\bf{X}}, {\bf{T}})}({\bf{x}}, {\bf{t}})$ denotes the marginal density function of $({\bf{X}},{\bf{T}})$. Therefore, the right-hand side of (\ref{A.1}) can be written as
		\begin{align*}
			&\frac{1}{2}n\det({\bf{H}}_n)[n_0({\bf{t}}_0){\bf{\Sigma}}_n({\bf{t}}_0)]^{-1/2}\int_{\mathbb{R}^p}\int_{\mathbb{R}^q}{\bf{x}}\sum_{l=1}^2\left\{{\bf{s}}^\top{\bf{H}}_n {\bm{\lambda}}^{(l)}({\bf{x}}, {\bf{t}}_0){\bf{H}}_n{\bf{s}}\right\}\\
			&\quad\times(1-F(\omega_n\mid {\bf{x}}, {\bf{t}}_0))f_{({\bf{X}}, {\bf{T}})}({\bf{x}}, {\bf{t}}_0)K({\bf{s}}){\rm{d}}{\bf{x}}{\rm{d}}{\bf{s}}\{1+o(1)\}\\
			&=[n_0({\bf{t}}_0){\bf{\Sigma}}_n({\bf{t}}_0)]^{1/2}\left\{\sum_{l=1}^2{\bf{\Lambda}}_n^{(l)}({\bf{t}}_0)\right\}\{1+o(1)\},
		\end{align*}
		where ${\bf{\Lambda}}_n^{(1)}({\bf{t}})$ and ${\bf{\Lambda}}_n^{(2)}({\bf{t}})$ are defined in Section \ref{subsec3.2}. Therefore, we have $E[[n_0({\bf{t}}_0){\bf{\Sigma}}_n({\bf{t}}_0)]^{-1/2}{\bf{A}}_n^{(2)}]$\\
		$-[n_0({\bf{t}}_0){\bf{\Sigma}}_n({\bf{t}}_0)]^{1/2}\sum_{l=1}^2{\bf{\Lambda}}_n^{(l)}({\bf{t}}_0)\to{\bf{0}}$ as $n\to\infty$. Hence, the proof of the first step is completed.
	\end{step}
	
	\begin{step}\label{step2}
		We abbreviate as
		\begin{equation*}
			{\bf{V}}_n({\bf{t}}_0)={\bf{X}}\left\{\gamma({\bf{X}}, {\bf{t}}_0)^{-1}\log\frac{Y}{\omega_n}-1\right\}I(Y>\omega_n)K({\bf{H}}_n^{-1}({\bf{t}}_0-{\bf{T}})).
		\end{equation*}
		Because $\{(Y_i, {\bf{X}}_i, {\bf{T}}_i)\}_{i=1}^n$ are independently distributed, we have
		\begin{align*}
			&{\rm{var}}[[n_0({\bf{t}}_0){\bf{\Sigma}}_n({\bf{t}}_0)]^{-1/2}\dot{\bm{L}}_n({\bm{\theta}}({\bf{t}}_0))]\\
			&=n[n_0({\bf{t}}_0){\bf{\Sigma}}_n({\bf{t}}_0)]^{-1/2}{\rm{var}}[{\bf{V}}_n({\bf{t}}_0)][n_0({\bf{t}}_0){\bf{\Sigma}}_n({\bf{t}}_0)]^{-1/2}\\
			&=n[n_0({\bf{t}}_0){\bf{\Sigma}}_n({\bf{t}}_0)]^{-1/2}E[{\bf{V}}_n({\bf{t}}_0){\bf{V}}_n({\bf{t}}_0)^\top][n_0({\bf{t}}_0){\bf{\Sigma}}_n({\bf{t}}_0)]^{-1/2}\\
			&\hspace*{0.5cm}-n[n_0({\bf{t}}_0){\bf{\Sigma}}_n({\bf{t}}_0)]^{-1/2}E[{\bf{V}}_n({\bf{t}}_0)]E[{\bf{V}}_n({\bf{t}}_0)]^\top[n_0({\bf{t}}_0){\bf{\Sigma}}_n({\bf{t}}_0)]^{-1/2}.
		\end{align*}
		From the result of Step \ref{step1}, the second term on the right-hand side converges to the zero matrix as $n\to\infty$. Using the Taylor expansion, the first term on the right-hand side can be written as
		\begin{align*}
			&n[n_0({\bf{t}}_0){\bf{\Sigma}}_n({\bf{t}}_0)]^{-1/2}E[{\bf{V}}_n({\bf{t}}_0){\bf{V}}_n({\bf{t}}_0)^\top][n_0({\bf{t}}_0){\bf{\Sigma}}_n({\bf{t}}_0)]^{-1/2}\\
			&=n[n_0({\bf{t}}_0){\bf{\Sigma}}_n({\bf{t}}_0)]^{-1/2}\\
			&\quad\times\left[\int_{\omega_n}^\infty\int_{\mathbb{R}^p}\int_{\mathbb{R}^q}{\bf{x}}{\bf{x}}^\top\left\{\gamma({\bf{x}}, {\bf{t}}_0)^{-1}\log\frac{y}{\omega_n}-1\right\}^2K({\bf{H}}_n^{-1}({\bf{t}}_0-{\bf{t}}))f(y, {\bf{x}}, {\bf{t}}){\rm{d}}y{\rm{d}}{\bf{x}}{\rm{d}}{\bf{t}}\right]\\
			&\quad\quad\times[n_0({\bf{t}}_0){\bf{\Sigma}}_n({\bf{t}}_0)]^{-1/2}\\
			&=n({\bf{t}}_0)[n_0({\bf{t}}_0){\bf{\Sigma}}_n({\bf{t}}_0)]^{-1/2}\\
			&\quad\times\left[\int_{\omega_n}^\infty\int_{\mathbb{R}^p}\int_{\mathbb{R}^q}{\bf{x}}{\bf{x}}^\top\left\{\gamma({\bf{x}}, {\bf{t}}_0)^{-1}\log\frac{y}{\omega_n}-1\right\}^2K({\bf{s}})\frac{f(y, {\bf{x}}, {\bf{t}}_0+{\bf{H}}_n{\bf{s}})}{f_{\bf{T}}({\bf{t}}_0)}{\rm{d}}y{\rm{d}}{\bf{x}}{\rm{d}}{\bf{s}}\right]\\
			&\quad\quad\times[n_0({\bf{t}}_0){\bf{\Sigma}}_n({\bf{t}}_0)]^{-1/2}\\
			&=\nu n({\bf{t}}_0)[n_0({\bf{t}}_0){\bf{\Sigma}}_n({\bf{t}}_0)]^{-1/2}\\
			&\quad\times\left[\int_{\omega_n}^\infty\int_{\mathbb{R}^p}{\bf{x}}{\bf{x}}^\top\left\{\gamma({\bf{x}}, {\bf{t}}_0)^{-1}\log\frac{y}{\omega_n}-1\right\}^2\frac{f(y, {\bf{x}}, {\bf{t}}_0)}{f_{\bf{T}}({\bf{t}}_0)}{\rm{d}}y{\rm{d}}{\bf{x}}\right]\\
			&\quad\quad\times[n_0({\bf{t}}_0){\bf{\Sigma}}_n({\bf{t}}_0)]^{-1/2}+o(1).\label{A.2}\tag{A.2}
		\end{align*}
		Because the conditional distribution of $\gamma({\bf{X}}, {\bf{t}}_0)^{-1}\log(Y/\omega_n)$ given $({\bf{X}}, {\bf{T}})=({\bf{x}}, {\bf{t}}_0)$ and $Y>\omega_n$ is approximately a standard exponential, we have
		\begin{equation*}
			\int_{\omega_n}^\infty\left\{\gamma({\bf{x}}, {\bf{t}}_0)^{-1}\log\frac{y}{\omega_n}-1\right\}^2\frac{f(y, {\bf{x}}, {\bf{t}}_0)/f_{({\bf{X}}, {\bf{T}})}({\bf{x}}, {\bf{t}}_0)}{1-F(\omega_n\mid {\bf{x}}, {\bf{t}}_0)}{\rm{d}}y\approx1.
		\end{equation*}
		Therefore, the right-hand side of (\ref{A.2}) can be written as
		\begin{align*}
			&\nu n({\bf{t}}_0)[n_0({\bf{t}}_0){\bf{\Sigma}}_n({\bf{t}}_0)]^{-1/2}\\
			&\quad\times\left[\int_{\mathbb{R}^p}{\bf{x}}{\bf{x}}^\top(1-F(\omega_n\mid {\bf{x}}, {\bf{t}}_0))\frac{f_{({\bf{X}}, {\bf{T}})}({\bf{x}}, {\bf{t}}_0)}{f_{\bf{T}}({\bf{t}}_0)}{\rm{d}}{\bf{x}}\right][n_0({\bf{t}}_0){\bf{\Sigma}}_n({\bf{t}}_0)]^{-1/2}+o(1)
		\end{align*}
		\begin{align*}
			&=\nu n({\bf{t}}_0)[n_0({\bf{t}}_0){\bf{\Sigma}}_n({\bf{t}}_0)]^{-1/2}\\
			&\quad\times\left[\int_{\mathbb{R}}\int_{\mathbb{R}^p}{\bf{x}}{\bf{x}}^\top I(y>\omega_n)\frac{f(y, {\bf{x}}, {\bf{t}}_0)}{f_{\bf{T}}({\bf{t}}_0)}{\rm{d}}y{\rm{d}}{\bf{x}}\right][n_0({\bf{t}}_0){\bf{\Sigma}}_n({\bf{t}}_0)]^{-1/2}+o(1)\\
			&\to\nu{\bf{I}}_p
		\end{align*}
		as $n\to\infty$. Hence, the proof of the second step is completed.
	\end{step}
	
	From the results of Steps \ref{step1} and \ref{step2}, we obtain Theorem \ref{thm1} by applying the Central Limit Theorem. 
\end{proof}

\begin{proof}[Proof of Theorem~{\upshape\ref{thm2}}]
	We show $E[[n_0({\bf{t}}_0){\bf{\Sigma}}_n({\bf{t}}_0)]^{-1/2}\ddot{\bm{L}}_n({\bm{\theta}}({\bf{t}}_0))[n_0({\bf{t}}_0){\bf{\Sigma}}_n({\bf{t}}_0)]^{-1/2}]\to{\bf{I}}_p$ as $n\to\infty$ in the third step and ${\rm{var}}[{\rm{vec}}([n_0({\bf{t}}_0){\bf{\Sigma}}_n({\bf{t}}_0)]^{-1/2}\ddot{\bm{L}}_n({\bm{\theta}}({\bf{t}}_0))[n_0({\bf{t}}_0){\bf{\Sigma}}_n({\bf{t}}_0)]^{-1/2})]\to{\bm{O}}$ as $n\to\infty$ in the fourth step, where ${\rm{vec}}(\cdot)$ is a vec operator.
	
	\begin{step}\label{step3}
		The conditional distribution of $\gamma({\bf{X}},{\bf{t}}_0)^{-1}\log(Y/\omega_n)$ given $({\bf{X}},{\bf{T}})=({\bf{x}},{\bf{t}}_0)$ and $Y>\omega_n$ is approximately a standard exponential. Consequently, we have
		\begin{align*}
			&E[[n_0({\bf{t}}_0){\bf{\Sigma}}_n({\bf{t}}_0)]^{-1/2}\ddot{\bm{L}}_n({\bm{\theta}}({\bf{t}}_0))[n_0({\bf{t}}_0){\bf{\Sigma}}_n({\bf{t}}_0)]^{-1/2}]\\
			&=E[{\bf{\Sigma}}_n({\bf{t}}_0)^{-1/2}\widetilde{\bf{\Sigma}}_n({\bf{t}}_0){\bf{\Sigma}}_n({\bf{t}}_0)^{-1/2}]+o(1).
		\end{align*}
		Under the condition (C.3), the right-hand side converges to ${\bf{I}}_p$ as $n\to\infty$. Hence, the proof of the third step is completed.
	\end{step}
	\begin{step}\label{step4}
		The conditional distribution of $\gamma({\bf{X}},{\bf{t}}_0)^{-1}\log(Y/\omega_n)$ given $({\bf{X}},{\bf{T}})=({\bf{x}},{\bf{t}}_0)$ and $Y>\omega_n$ is approximately a standard exponential. Consequently, we have
		\begin{align*}
			&{\rm{var}}[{\rm{vec}}([n_0({\bf{t}}_0){\bf{\Sigma}}_n({\bf{t}}_0)]^{-1/2}\ddot{\bm{L}}_n({\bm{\theta}}({\bf{t}}_0))[n_0({\bf{t}}_0){\bf{\Sigma}}_n({\bf{t}}_0)]^{-1/2})]\\
			&={\rm{var}}[({\bf{\Sigma}}_n({\bf{t}}_0)^{-1/2}\otimes{\bf{\Sigma}}_n({\bf{t}}_0)^{-1/2}){\rm{vec}}(n_0({\bf{t}}_0)^{-1/2}\ddot{\bm{L}}_n({\bm{\theta}}({\bf{t}}_0)))]\\
			&=({\bf{\Sigma}}_n({\bf{t}}_0)^{-1/2}\otimes{\bf{\Sigma}}_n({\bf{t}}_0)^{-1/2})\\
			&\quad\times{\rm{var}}[{\rm{vec}}(n_0({\bf{t}}_0)^{-1/2}\ddot{\bm{L}}_n({\bm{\theta}}({\bf{t}}_0)))]({\bf{\Sigma}}_n({\bf{t}}_0)^{-1/2}\otimes{\bf{\Sigma}}_n({\bf{t}}_0)^{-1/2})\\
			&=({\bf{\Sigma}}_n({\bf{t}}_0)^{-1/2}\otimes{\bf{\Sigma}}_n({\bf{t}}_0)^{-1/2})\\
			&\quad\times{\rm{var}}[{\rm{vec}}(\widetilde{\bf{\Sigma}}_n({\bf{t}}_0))]({\bf{\Sigma}}_n({\bf{t}}_0)^{-1/2}\otimes{\bf{\Sigma}}_n({\bf{t}}_0)^{-1/2})+o(1)\\
			&={\rm{var}}[({\bf{\Sigma}}_n({\bf{t}}_0)^{-1/2}\otimes{\bf{\Sigma}}_n({\bf{t}}_0)^{-1/2}){\rm{vec}}(\widetilde{\bf{\Sigma}}_n({\bf{t}}_0))]+o(1)\\
			&={\rm{var}}[{\rm{vec}}({\bf{\Sigma}}_n({\bf{t}}_0)^{-1/2}\widetilde{\bf{\Sigma}}_n({\bf{t}}_0){\bf{\Sigma}}_n({\bf{t}}_0)^{-1/2})]+o(1),
		\end{align*}
		where ``$\otimes$'' denotes the Kronecker product. Under the condition (C.3), the right-hand side converges to ${\bm{O}}$ as $n\to\infty$. Hence, the proof of the fourth step is completed.
	\end{step}
\end{proof}

\begin{proof}[Proof of Theorem~{\upshape\ref{thm3}}]
	We define ${\bm{\alpha}}_n={\bf{\Sigma}}_n({\bf{t}}_0)^{1/2}({\bm{\theta}}-{\bm{\theta}}({\bf{t}}_0))$ and ${\bm{\alpha}}_n^*={\bf{\Sigma}}_n({\bf{t}}_0)^{1/2}{\bm{\theta}}({\bf{t}}_0)$. Additionally, we define ${\bf{Z}}_{ni}={\bf{\Sigma}}_n({\bf{t}}_0)^{-1/2}{\bf{X}}_i,\ i=1, \ldots, n$. The objective function $L_n({\bm{\theta}})$ can be written as
	\begin{align*}
		L_n^*({\bm{\alpha}}_n)&=\sum_{i=1}^n\left\{\exp[{\bf{Z}}_{ni}^\top({\bm{\alpha}}_n+{\bm{\alpha}}_n^*)]\log\frac{Y_i}{\omega_n}-{\bf{Z}}_{ni}^\top({\bm{\alpha}}_n+{\bm{\alpha}}_n^*)\right\}\\
		&\quad\times I(Y_i>\omega_n)K({\bf{H}}_n^{-1}({\bf{t}}_0-{\bf{T}}_i)).
	\end{align*}
	Let $\dot{\bm{L}}_n^*({\bm{\alpha}}_n)$ be the gradient vector of $L_n^*({\bm{\alpha}}_n)$ and $\ddot{\bm{L}}_n^*({\bm{\alpha}}_n)$ be Hessian matrix of $L_n^*({\bm{\alpha}}_n)$. We assume $\ddot{\bm{L}}_n^*({\bm{\alpha}}_n)$ is a positive definite matrix for all ${\bm{\alpha}}_n\in\mathbb{R}^p$. Therefore, $L_n^*({\bm{\alpha}}_n)$ is strictly convex.
	
	Using the Taylor expansion, we have
	\begin{align*}
		L_n^*\left(\frac{{\bf{s}}}{\sqrt{n_0({\bf{t}}_0)}}\right)=L_n^*({\bm{0}})+{\bf{s}}^\top\frac{\dot{\bm{L}}_n^*({\bm{0}})}{\sqrt{n_0({\bf{t}}_0)}}+\frac{1}{2}{\bf{s}}^\top\frac{\ddot{\bm{L}}_n^*({\bm{0}})}{n_0({\bf{t}}_0)}{\bf{s}}+o_P(1)
	\end{align*}
	for fixed ${\bf{s}}\in\mathbb{R}^p$. By applying Theorems \ref{thm1} and \ref{thm2}, we have $\dot{\bm{L}}_n^*({\bm{0}})/\sqrt{n_0({\bf{t}}_0)}+[n_0({\bf{t}}_0){\bm{\Sigma}}_n({\bf{t}}_0)]^{1/2}\sum_{l=1}^2{\bf{\Lambda}}_n^{(l)}({\bf{t}}_0)\xrightarrow{D}{\rm{N}}(-{\bf{b}}({\bf{t}}_0),\nu{\bf{I}}_p)$ and $\ddot{\bm{L}}_n^*({\bm{0}})/n_0({\bf{t}}_0)\xrightarrow{P}{\bf{I}}_p$. We assume $[n_0({\bf{t}}_0){\bm{\Sigma}}_n({\bf{t}}_0)]^{1/2}\sum_{l=1}^2{\bf{\Lambda}}_n^{(l)}({\bf{t}}_0)\to0$ as $n\to\infty$. Consequently, this implies that, for any $\varepsilon>0$, there exists a constant $C$ such that
	\begin{equation*}
		\liminf_{n\to\infty}P\left(\inf_{{\bf{s}}\in\mathbb{R}^p:\|{\bf{s}}\|=C}L_n^*\left(\frac{{\bf{s}}}{\sqrt{n_0({\bf{t}}_0)}}\right)>L_n^*({\bm{0}})\right)>1-\varepsilon,
	\end{equation*}
	which implies $\widehat{\bm{\alpha}}_n={\bf{\Sigma}}_n({\bf{t}}_0)^{1/2}(\widehat{\bm{\theta}}({\bf{t}}_0)-{\bm{\theta}}({\bf{t}}_0))=O_P(\sqrt{n_0({\bf{t}}_0)})$ (Fan and Li 2001). From the Taylor expansion of $\dot{\bm{L}}_n^*(\widehat{\bm{\alpha}}_n)={\bm{0}}$, we have
	\begin{equation*}
		{\bm{0}}=\frac{\dot{\bm{L}}_n^*\left({\bm{0}}\right)}{\sqrt{n_0({\bf{t}}_0)}}+\sqrt{n_0({\bf{t}}_0)}\widehat{\bm{\alpha}}_n^\top\frac{\ddot{\bm{L}}_n^*\left({\bm{0}}\right)}{n_0({\bf{t}}_0)}+o_P\left(1\right).
	\end{equation*}
	Therefore, by applying Theorems \ref{thm1} and \ref{thm2}, we obtain Theorem \ref{thm3}.
\end{proof}


\begin{thebibliography}{99}
	\bibitem{}Andriyana, Y., Gijbels, I., and Verhasselt, A. (2014). ``P-splines quantile regression estimation in varying coefficient models''. {\it{TEST}} {\bf{23}}, 153-194. \url{https://doi.org/10.1007/s11749-013-0346-2}
	\bibitem{}Andriyana, Y., Gijbels, I., and Verhasselt, A. (2018). ``Quantile regression in varying-coefficient models: non-crossing quantile curves and heteroscedasticity''. {\it{Statistical Papers}} {\bf{59}}, 1589-1621. \url{https://doi.org/10.1007/s00362-016-0847-7}
	\bibitem{}Cai, Z., and Xu, X. (2008). ``Nonparametric Quantile Estimations for Dynamic Smooth Coefficient Models''. {\it{Journal of the American Statistical Association}} {\bf{103}}, 1595-1608. \url{http://dx.doi.org/10.1198/016214508000000977}
	\bibitem{}Cai, Z., Fan, J., and Yao, Q. (2000). ``Functional-Coefficient Regression Models for Nonlinear Time Series''. {\it{Journal of the American Statistical Association}} {\bf{95}}, 941-956. \url{http://dx.doi.org/10.1080/01621459.2000.10474284}
	\bibitem{}Csorgo, S., Deheuvels, P., and Mason, D. (1985). ``Kernel Estimates of the Tail Index of a Distribution''. {\it{The Annuals of Statistics}} {\bf{13}}, 1050-1077. \url{http://dx.doi.org/10.1214/aos/1176349656}
	\bibitem{}de Haan, L., and Ferreira, A. (2006). {\it{Extreme Value Theory: An Introduction}}. Springer, New York. \url{https://doi.org/10.1007/0-387-34471-3}
	\bibitem{}de Haan, L., and Zhou, C. (2022). ``Bootstrapping Extreme Value Estimators'', {\it{Journal of the American Statistical Association}}, {\bf{0}}, 1-12. \url{https://doi.org/10.1080/01621459.2022.2120400}
	\bibitem{}Daouia, A., Gardes, L., and Girard, S. (2013). ``On kernel smoothing for extremal quantile regression''. {\it{Bernoulli}} {\bf{19}}, 2557-2589. \url{http://dx.doi.org/10.3150/12-BEJ466}
	\bibitem{}Dekkers, A. L. M., Einmahl, J. H. J., and de Haan, L. (1989). ``A Moment Estimator for the Index of an Extreme-Value Distribution''. {\it{The Annals of Statistics}} {\bf{17}}, 1833-1855. \url{http://dx.doi.org/10.1214/aos/1176347397}
	\bibitem{}Fan, J., and Li, R. (2001). ``Variable Selection via Nonconcave Penalized Likelihood and its Oracle Properties''. {\it{Journal of the American Statistical Association}} {\bf{96}}, 1348-1360. \url{http://dx.doi.org/10.1198/016214501753382273}
	\bibitem{}Fan, J., and Zhang, W. (1999). ``Statistical estimation in varying coefficient models''. {\it{The Annals of Statistics}} {\bf{27}}, 1491-1518. \url{http://dx.doi.org/10.1214/aos/1017939139}
	\bibitem{}Fan, J., and Zhang, W. (2000). ``Simultaneous Confidence Bands and Hypothesis Testing in Varying-coefficient Models''. {\it{Scandinavian Journal of Statistics}} {\bf{27}}, 715-731. \url{http://dx.doi.org/10.1111/1467-9469.00218}
	\bibitem{}Fan, J., and Zhang, W. (2008). ``Statistical methods with varying coefficient models''. {\it{Statistics and Its Interface}} {\bf{1}}, 179-195. \url{http://dx.doi.org/10.4310/SII.2008.v1.n1.a15}
	\bibitem{}Fan, J., Zhang, C., and Zhang, J. (2001). ``Generalized Likelihood Ratio Statistics and Wilks Phenomenon''. {\it{The Annals of Statistics}} {\bf{29}}, 153-193. \url{http://dx.doi.org/10.1214/aos/996986505}
	\bibitem{}Gardes, L., and Girard, S. (2010). ``Conditional extremes from heavy-tailed distributions: Application to the estimation of extreme rainfall return levels''. {\it{Extremes}} {\bf{13}}, 177-204. \url{http://dx.doi.org/10.1007/s10687-010-0100-z}
	\bibitem{}Gardes, L., and Stupfler, G. (2014). ``Estimation of the conditional tail index using a smoothed local hill estimator''. {\it{Extremes}} {\bf{17}}, 45-75. \url{http://dx.doi.org/10.1007/s10687-013-0174-5}
	\bibitem{}Goegebeur, Y., Guillou, A., and Schorgen, A. (2014). ``Nonparametric regression estimation of conditional tails: random covariate case''. {\it{Statistics}} {\bf{48}}, 732-755. \url{https://doi.org/10.1080/02331888.2013.800064}
	\bibitem{}Goegebeur, Y., Guillou, A., and Stupfler, G. (2015). ``Uniform asymptotic properties of the nonparametric regression estimator of conditional tails''. {\it{Annales de l' Institut Henri Poincar\'{e}, Probabilit\'{e}s et Statistiques}} {\bf{51}}, 1190-1213. \url{http://dx.doi.org/10.1214/14-AIHP624}
	\bibitem{}Gomes, M., de Haan., L., and Peng., L (2002). ``Semi-parametric Estimation of the Second Order Parameter in Statistics of Extremes'', {\it{Extremes}} {\bf{5}}, 387-414. \url{https://doi.org/10.1023/A:1025128326588}
	\bibitem{}Hall, P. (1982). ``On Some Simple Estimates of an Exponent of Regular Variation''. {\it{Journal of the Royal Statistical Society: Series B (Methodological)}} {\bf{44}}, 37-42. \url{http://dx.doi.org/10.1111/j.2517-6161.1982.tb01183.x}
	\bibitem{}Hill, B. M. (1975). ``A Simple General Approach to Inference About the Tail of a Distribution''. {\it{The Annuals of Statistics}} {\bf{3}}, 1163-1174. \url{http://dx.doi.org/10.1214/aos/1176343247}
	\bibitem{}Hastie, T., and Tibshirani, R. (1993). ``Varying-Coefficient Models''. {\it{Journal of the Royal Statistical Society: Series B (Methodological)}} {\bf{55}}, 757-779. \url{http://dx.doi.org/10.1111/j.2517-6161.1993.tb01939.x}
	\bibitem{}Huang, J. Z., Wu, C. O., and Zhou, L. (2002). ``Varying-coefficient models and basis function approximations for the analysis of repeated measurements''. {\it{Biometrika}} {\bf{89}}, 111-128. \url{http://dx.doi.org/10.1093/biomet/89.1.111}
	\bibitem{}Huang, J. Z., Wu, C. O., and Zhou, L. (2004). ``Polynomial Spline Estimation and Inference for Varying Coefficient Models with Longitudinal Data''. {\it{Statistica Sinica}} {\bf{14}}, 763-788.
	\bibitem{}Kim, M. O. (2007). ``Quantile regression with varying coefficients''. {\it{The Annals of Statistics}} {\bf{35}}, 92-108. \url{http://dx.doi.org/10.1214/009053606000000966}
	\bibitem{}Li, R., Leng, C., and You, J. (2022). ``Semiparametric Tail Index Regression''. {\it{Journal of Business \& Economic Statistics}} {\bf{40}}, 82-95. \url{http://dx.doi.org/10.1080/07350015.2020.1775616}
	\bibitem{}Ma, Y., Jiang, Y., and Huang, W. (2019). ``Tail index varying coefficient model''. {\it{Communications in Statistics - Theory and Methods}} {\bf{48}}, 235-256. \url{http://dx.doi.org/10.1080/03610926.2017.1406519}
	\bibitem{}Ma, Y., Wei, B., and Huang, W. (2020). ``A nonparametric estimator for the conditional tail index of Pareto-type distributions''. {\it{Metrika}} {\bf{83}}, 17-44. \url{http://dx.doi.org/10.1007/s00184-019-00723-8}
	\bibitem{}Pickands, J. (1975). ``Statistical Inference Using Extreme Order Statistics''. {\it{The Annals of Statistics}} {\bf{3}}, 119-131. \url{http://dx.doi.org/10.1214/aos/1176343003}
	\bibitem{}Rosenblatt, M. (1976). ``On the Maximal Deviation of $k$-Dimensional Density Estimates''. {\it{The Annals of Probability}} {\bf{4}}, 1009-1015. \url{http://dx.doi.org/10.1214/aop/1176995945}
	\bibitem{}Smith, R. L. (1987). ``Estimating Tails of Probability Distributions''. {\it{The Annals of Statistics}} {\bf{15}}, 1174-1207. \url{http://dx.doi.org/10.1214/aos/1176350499}
	\bibitem{}Stupfler, G. (2013). ``A moment estimator for the conditional extreme-value index''. {\it{Electronic Journal of Statistics}} {\bf{7}}, 2298-2343. \url{http://dx.doi.org/10.1214/13-EJS846}
	\bibitem{}Wang, H., and Tsai, C. L. (2009). ``Tail Index Regression''. {\it{Journal of the American Statistical Association}} {\bf{104}}, 1233-1240. \url{http://dx.doi.org/10.1198/jasa.2009.tm08458}
	\bibitem{}Wu, C. O., Chiang, C. T., and Hoover, D. R. (1998). ``Asymptotic Confidence Regions for Kernel Smoothing of a Varying-Coefficient Model with Longitudinal Data''. {\it{Journal of the American Statistical Association}} {\bf{93}}, 1388-1402. \url{http://dx.doi.org/10.1080/01621459.1998.10473800}
	\bibitem{}Yoshida, T. (2023). ``Single-index models for extreme value index regression''. {\it arXiv} \url{http://arxiv.org/abs/2203.05758}
	\bibitem{}Youngman, B. D. (2019). ``Generalized Additive Models for Exceedances of High Thresholds With an Application to Return Level Estimation for U.S. Wind Gusts''. {\it{Journal of the American Statistical Association}} {\bf{114}}, 1865-1879. \url{http://dx.doi.org/10.1080/01621459.2018.1529596}
\end{thebibliography}
\end{document}